\documentclass[12pt, reqno]{amsart}
\usepackage{amsmath, amstext, amsbsy, amssymb, amscd}

\newtheorem{theorem}{Theorem}[section]
\newtheorem{lemma}[theorem]{Lemma}
\newtheorem{proposition}[theorem]{Proposition}
\newtheorem{corollary}[theorem]{Corollary}
\theoremstyle{definition}
\newtheorem{definition}[theorem]{Definition}
\newtheorem{example}[theorem]{Example}

\theoremstyle{remark}
\newtheorem{remark}[theorem]{Remark}

\numberwithin{equation}{section}

\newcommand{\ch}{\mbox{ch} }
\newcommand{\C}{ \mathbb C }

\newcommand{\End}{{\rm End}}
\newcommand{\fock}{{\mathbb H}_X}
\newcommand{\g}{{\gamma}}
\newcommand{\In}{\mathcal I^{[n]}}

\newcommand{\KX}{K}
\newcommand{\lambsq}{s(\lambda)}
\newcommand{\orbsymbar}{H^*_{\text{orb}}(\overline{X}^n/S_n)}

\newcommand{\orbsym}{H^*_{\text{orb}}(X^n/S_n)}
\newcommand{\Hn}{H^*(\Xn)}

\newcommand{\vac}{|0\rangle}

\newcommand{\Xbar}{ \overline{X}}
\newcommand{\Xn}{ X^{[n]}}
\newcommand{\Z}{ \mathbb Z }

\newenvironment{demo}[1]%
{\vskip-\lastskip\medskip
  \noindent
  {\em #1.}\enspace
  }%
{\qed\par\medskip
  }

\begin{document}
\title[Ideals of the cohomology rings and applications]
    {Ideals of the cohomology rings of Hilbert schemes
and their applications}
\author[Wei-Ping Li]{Wei-Ping Li$^1$}
\address{Department of Mathematics, HKUST, Clear Water Bay, Kowloon, Hong
Kong } \email{mawpli@ust.hk}
\thanks{${}^1$Partially supported by the grant HKUST6170/99P}

\author[Zhenbo Qin]{Zhenbo Qin$^2$}
\address{Department of Mathematics, University of Missouri, Columbia, MO
65211, USA} \email{zq@math.missouri.edu}
\thanks{${}^2$Partially supported by an NSF grant}

\author[Weiqiang Wang]{Weiqiang Wang$^3$}
\address{Department of Mathematics, University of Virginia,
Charlottesville, VA 22904} \email{ww9c@virginia.edu}
\thanks{${}^3$Partially supported by an NSF grant}


\begin{abstract}
We study the ideals of the rational cohomology ring of the Hilbert
scheme $\Xn$ of $n$ points on a smooth projective surface $X$.
%which are induced from ideals of the cohomology ring of $X$.
%We prove that part of the cup product structure constants
%for the cohomology rings of $\Xn$ with respect to
%the Heisenberg monomial linear basis is independent of $n$.
As an application, for a large class of smooth quasi-projective
surfaces $X$, we show that every cup product structure constant of
$H^*(\Xn)$ is independent of $n$; moreover, we obtain two sets of
ring generators for the cohomology ring $H^*(\Xn)$.

Similar results are established for the Chen-Ruan orbifold
cohomology ring of the symmetric product. In particular, we prove
a ring isomorphism between $H^*(\Xn; \C)$ and $H^*_{\rm
orb}(X^n/S_n; \C)$ for a large class of smooth quasi-projective
surfaces with numerically trivial canonical class.
\end{abstract}

\maketitle
\date{}
%\tableofcontents

%%
%%
%%
%%
%%
%%
%%
\section{Introduction}

This is a sequel to [LQW1-4] and \cite{QW}. We continue the study
of the cohomology rings of the Hilbert schemes $\Xn$ of $n$ points
on a smooth surface $X$ and the Chen-Ruan orbifold cohomology
rings of the symmetric products $X^n/S_n$. In this paper, the
surface $X$ is allowed to be projective as well as
quasi-projective (our usage of the terminology ``quasi-projective"
excludes ``projective").

In \cite{Lehn} [LQW1-4] \cite{LS2} which were in turn built on the
earlier works \cite{Got, VW, Na1, Na2, Gro} and others, the
connections between vertex operators and the multiplicative
structure of the rational cohomology group $H^*(\Xn)$ when $X$ is
projective have been developed. These connections have been
successfully applied to unravel various structures on the
cohomology ring of $\Xn$. However, the situation changes
dramatically when $X$ is quasi-projective. To date, the
understanding of the cohomology ring $H^*(\Xn)$ for $X$
quasi-projective has been rather limited with the exception of the
affine plane \cite{ES, Lehn, LS1, Vas} and the minimal resolution
of a simple singularity $\mathbb C^2/\Gamma$ where $\Gamma$ is a
finite group of $SL_2(\mathbb C)$ (cf. \cite{Wa}).

Besides the minimal resolutions just mentioned above, typical
important examples of quasi-projective surfaces include the
cotangent bundle of a smooth projective curve and the surface
obtained from a smooth projective surface by deleting a point. All
these surfaces are among a class of quasi-projective surfaces
which satisfy what we call the {\em S-property} (see
Definition~\ref{S}). One of the goals of the present paper is to
establish some  general results of the cohomology rings of the
Hilbert schemes of points for such a large class of
quasi-projective surfaces. These results are in general not valid
for projective surfaces, and conjecturally, they hold for every
quasi-projective surface (without the S-property assumption).

We begin with a study of the ideals $\In$ of the cohomology ring
$H^*(\Xn)$ for a smooth {\em projective} surface $X$, which are
induced from ideals $\mathcal I$ of the cohomology ring of $X$
(see Definition~\ref{pushforward_ideal}). We prove that part of
the cup product structure constants for the cohomology rings of
the Hilbert schemes $\Xn$ with respect to the
Nakajima-Grojnowski's Heisenberg monomial linear basis is
independent of $n$. This part of the cup product structure
constants can be regarded as coming from the complement to the
ideal $\In$. The methods used in establishing the above results
follow the techniques developed in \cite{LQW3, LQW4}.

While the results on the ideals $\In$ are of independent interest,
they are developed with applications to smooth quasi-projective
surfaces in mind. As the first application, we prove the following
result (Theorem~\ref{structure_constant_Y}).

\begin{theorem} \label{intro_thm1}
Let $X$ be a smooth quasi-projective surface with the S-property.
Then, all the cup product structure constants of $H^*(\Xn)$ are
independent of $n$.
\end{theorem}

The precise definition of the structure constants is given by
(\ref{str_constant_Y}). The $n$-independence has been conjectured
in \cite{Wa} to be true for every smooth quasi-projective surface,
and has its counterpart in the framework of the class algebras of
the wreath products. It allows us to introduce a universal ring
$\mathcal G_X$ for a smooth quasi-projective surface $X$ with the
S-property, which governs the cohomology rings $H^*(\Xn)$ for a
fixed $X$ and all $n$. This ring $\mathcal G_X$ will be called
{\it the FH ring associated to $X$} as it is analogous to the one
arising in the framework of symmetric groups and wreath products
\cite{FH, Wa}. We further determine the ring structure of
$\mathcal G_X$, and obtain two sets of ring generators of
$H^*(\Xn)$ which are the quasi-projective counterparts of the main
results in \cite{LQW1, LQW2}. We remark that a universal ring
termed as the Hilbert ring was introduced in \cite{LQW3} which
governs the cohomology rings $H^*(\Xn)$ for a fixed {\it
projective} $X$ and all $n$. These two universal rings reflect
distinct structures of the corresponding cohomology rings.

Observe that a distinguished Heisenberg generator (i.e. the first
annihilation operator) $\frak a_1([x]_c)$ is cohomology
degree-preserving, where $[x]_c \in H^4_c(X)$ denotes the
Poincar\'e dual of the homology class in $H_0(X)$ represented by a
point in $X$. It turns out that Theorem~\ref{intro_thm1} admits
the following equivalent reformulation.

\begin{corollary}
Let $X$ be a smooth quasi-projective surface with the S-property.
Then $-\frak a_1([x]_c)\colon  H^*(X^{[n+1]}) \rightarrow
H^*(X^{[n]})$ is a surjective ring homomorphism.
\end{corollary}

Recall a well-known fact that the Hilbert-Chow morphism from the
Hilbert schemes $\Xn$ to the symmetric product $X^n/S_n$ is a
resolution of singularity. In \cite{CR}, Chen and Ruan introduced
the notion of the orbifold cohomology ring $H_{\text{orb}}^*(Y)$
for an orbifold $Y$. Following the general machinery developed in
\cite{QW}, we automatically establish results parallel to those
stated in the previous paragraphs for the orbifold cohomology
rings $H^*_{\rm orb}(X^n/S_n)$ of the symmetric products
$X^n/S_n$. We further prove the following result
(Theorem~\ref{th:ringisom} and Remark~\ref{rmk:ringisom}).

\begin{theorem} \label{intro_thm2}
Let $X$ be a smooth quasi-projective surface with the S-property
and a numerically trivial canonical class. Then, the cohomology
ring $H^*(\Xn; \C)$ is isomorphic to the orbifold cohomology ring
$H^*_{\text{\rm orb}}(X^n/S_n; \C)$.
\end{theorem}

A conjecture of Ruan \cite{Ru1, Ru2} states that the cohomology
ring $H^*(Z; \mathbb C)$ with $\C$-coefficient is isomorphic to
the orbifold cohomology ring $H_{\text{orb}}^*(Y; \mathbb C)$ with
$\C$-coefficient for any hyperkahler resolution $Z$ of an orbifold
$Y$. In light of this conjecture, it is a very interesting
question for which surfaces satisfying the assumption in
Theorem~\ref{intro_thm2} the corresponding Hilbert schemes of
points carry a hyperkahler structure, (and
Theorem~\ref{intro_thm2} for these surfaces confirms Ruan's
conjecture). For example, the Hilbert scheme of points on the
minimal resolution of a simple singularity carries a hyperkahler
structure \cite{Na2}. We remark that Theorem~\ref{intro_thm2}
holds when the surface $X$ is the cotangent bundle of a smooth
projective curve, or the minimal resolution of a simple
singularity $\mathbb C^2/\Gamma$. In the special case when $X$ is
the affine plane (i.e. $\Gamma$ is trivial), we recover the main
results of \cite{LS1, Vas} (also cf. \cite{Ru1}). An isomorphism
as in Theorem~\ref{intro_thm2} when $X$ is projective with a
numerically trivial canonical class was earlier established by
various authors \cite{LS2, FG, Uri, QW}, which also supported
Ruan's conjecture.

The paper is organized as follows. In Sect.~\ref{sect_ideals}, we
study the ideals in $H^*(\Xn)$ for a smooth projective surface
$X$. In Sect.~\ref{sect_independence}, we verify the partial
$n$-independence of the cup product structure constants for the
cohomology rings of $\Xn$ when $X$ is a smooth projective surface.
In Sect.~\ref{sect_quasi}, we prove Theorem~\ref{intro_thm1}, and
construct as well as study the FH ring $\mathcal G_X$. In
Sect.~\ref{sect_orb}, for a smooth quasi-projective surface $X$
with the S-property, we formulate the analogous results for the
orbifold cohomology ring of the symmetric product $X^n/S_n$.
Furthermore, we establish Theorem~\ref{intro_thm2}.

\medskip\noindent
{\bf Conventions:} All cohomology groups are in $\mathbb
Q$-coefficients unless otherwise indicated. The cup product of two
cohomology classes $\alpha$ and $\beta$ is denoted by $\alpha
\cdot \beta$ or simply by $\alpha \beta$. For a continuous map $p:
Y_1 \to Y_2$ between two smooth compact manifolds and for
$\alpha_1 \in H^*(Y_1)$, let $p_*(\alpha_1) =
\text{PD}^{-1}p_{*}(\text{PD}(\alpha_1))$ where $\text{PD}$ stands
for the Poincar\'e duality. For a smooth projective surface $X$,
let $1_X \in H^0(X)$ be the fundamental cohomology class of $X$,
and $[x] \in H^4(X)$ be the Poincar\'e dual of the homology class
represented by a point $x\in X$.

\medskip\noindent
{\bf Acknowledgment.} We thank M. Lehn and Y. Ruan for helpful
discussions.

\section{Ideals in $H^*(\Xn)$ for $X$ projective}
\label{sect_ideals}

\subsection{Preliminaries}
\par
${}$

Let $X$ be a smooth projective complex surface with the canonical
class $\KX$ and the Euler class $e$, and $\Xn$ be the Hilbert
scheme of points in $X$. It is known that $X^{[n]}$ is a
desingularization of the symmetric product $X^n/S_n$. Let
$\mathcal Z_n\subset X^{[n]}\times X$ be the universal subscheme,
and $\fock = \oplus_{n=0}^\infty \Hn$. Recall that a Heisenberg
algebra was defined in \cite{Na2} acting on $\fock$. In this paper
we follow the notations and conventions in \cite{LQW3, LQW4} for
the Heisenberg algebra generators however. Namely, we have the
operators $\mathfrak a_{n}(\alpha) \in \End(\fock)$ with $\alpha
\in H^*(X)$ and $n \in \mathbb Z$ which satisfy the following
Heisenberg algebra commutation relation:
\begin{eqnarray}  \label{eq:heis}
[\mathfrak a_m(\alpha), \mathfrak a_n(\beta)] = -m \;
\delta_{m,-n} \int_X(\alpha \beta) \cdot {\rm Id}_{\fock}.
\end{eqnarray}
Here and throughout the paper, the Lie brackets are understood in
the super sense according to the parity of the cohomology degrees
of the cohomology classes involved.
%Put $\mathfrak a_0(\alpha) =0$.
When $n > 0$, we often refer to $\mathfrak a_{-n}(\alpha)$ (resp.
$\mathfrak a_n(\alpha)$) as the {\em creation} (resp. {\em
annihilation}) operator. The space $\fock$ is an irreducible
module over the Heisenberg algebra generated by the operators
$\mathfrak a_n(\alpha)$ with a highest~weight~vector $\vac=1 \in
H^0(X^{[0]}) \cong \mathbb Q$. It follows that $\fock$ is linearly
spanned by all the Heisenberg monomials $\mathfrak
a_{-n_1}(\alpha_1) \cdots \mathfrak a_{-n_k}(\alpha_k)\vac$ where
$k \ge 0$ and $n_1, \ldots, n_k > 0$.

For $n > 0$ and a homogeneous class $\gamma \in H^*(X)$, let
$|\gamma| = s$ if $\gamma \in H^s(X)$, and let $G_i(\gamma, n)$ be
the homogeneous component in $H^{|\gamma|+2i}(\Xn)$ of
\begin{eqnarray*}
G(\gamma, n) = p_{1*}(\ch({\mathcal O}_{{\mathcal Z}_n}) \cdot
p_2^*{\rm td}(X) \cdot p_2^*\gamma) \in \Hn
\end{eqnarray*}
where $p_1$ and $p_2$ be the projections of $\Xn \times X$ to
$\Xn$ and $X$ respectively. We extend the notion $G_i(\gamma, n)$
linearly to an arbitrary $\gamma \in H^*(X)$, and set $G(\gamma,
0) =0$. The {\it Chern character operator} ${\mathfrak
G}_i(\gamma) \in \End({\fock})$ is defined to be the operator
acting on the component $H^*(\Xn)$ by the cup product with
$G_i(\gamma, n)$. It was proved in \cite{LQW1} that the cohomology
ring of $\Xn$ is generated by the classes $G_{i}(\gamma, n)$ where
$0 \le i < n$ and $\gamma$ runs over a linear basis of $H^*(X)$.

For $k \ge 1$, let $\tau_{k*}: H^*(X) \to H^*(X^k)$ be the map
induced by the diagonal embedding $\tau_k: X \to X^k$, and let
$\mathfrak a_{m_1} \cdots \mathfrak a_{m_k}(\tau_{k*}(\alpha))$
denote $\sum_j \mathfrak a_{m_1}(\alpha_{j,1}) \cdots \mathfrak
a_{m_k}(\alpha_{j,k})$ when $\tau_{k*}\alpha = \sum_j \alpha_{j,1}
\otimes \cdots \otimes \alpha_{j, k}$ via the K\"unneth
decomposition of $H^*(X^k)$.

The following two lemmas were proved in \cite{LQW3}, where
$\tau_{0*}(\alpha)$ denotes $\int_X \alpha$.

\begin{lemma} \label{k_s}
Let $k, s \ge 1$, $n_1, \ldots, n_k, m_1, \ldots, m_s \in \Z$, and
$\alpha, \beta \in H^*(X)$.
 \begin{enumerate}
\item[{\rm (i)}] The commutator
$[\mathfrak a_{n_1} \cdots \mathfrak a_{n_{k}} (\tau_{k*}\alpha),
\mathfrak a_{m_1} \cdots \mathfrak a_{m_{s}}(\tau_{s*}\beta)]$ is
equal to
\begin{eqnarray*}
-\sum_{t=1}^k \sum_{j=1}^s n_t \delta_{n_t,-m_j} \cdot \left (
\prod_{\ell=1}^{j-1} \mathfrak a_{m_\ell} \prod_{1 \le u \le k, u
\ne t} \mathfrak a_{n_u} \prod_{\ell=j+1}^{s} \mathfrak a_{m_\ell}
\right )(\tau_{(k+s-2)*}(\alpha\beta));
\end{eqnarray*}

\item[{\rm (ii)}] Let $j$ satisfy $1 \le j < k$. Then,
$\mathfrak a_{n_1} \cdots \mathfrak a_{n_k}(\tau_{k*}\alpha)$ is
equal to
\begin{eqnarray*}
\left ( \prod_{1 \le s < j} \mathfrak a_{n_s} \cdot \mathfrak
a_{n_{j+1}} \mathfrak a_{n_{j}} \cdot \prod_{j+1 < s \le k}
\mathfrak a_{n_s} \right ) (\tau_{k*}\alpha) - n_j
\delta_{n_j,-n_{j+1}} \prod_{1 \le s \le k \atop s \ne j, j+1}
\mathfrak a_{n_s}(\tau_{(k-2)*}(e\alpha)).
\end{eqnarray*}
 \end{enumerate}
\end{lemma}

\begin{lemma} \label{nonsense1}
Fix $k \ge 0$ and $b \ge 1$. Let $\frak g \in \End(\mathbb H)$ be
of bi-degree $(\tilde s, s)$ satisfying
\begin{eqnarray} \label{nonsense1.1}
[[\cdots [\frak g, \frak a_{m_1}(\beta_1)], \cdots], \frak
a_{m_{k+2}}(\beta_{k+2})] =0
\end{eqnarray}
whenever $m_i < 0$ and $\beta_i \in H^*(X)$ for each $i$. Let $A =
\mathfrak a_{-n_1}(\alpha_{1}) \cdots \mathfrak
a_{-n_b}(\alpha_{b})|0\rangle$ where $n_1, \ldots, n_b > 0$ and
$\alpha_{1}, \ldots, \alpha_b \in H^*(X)$. Then, $\frak g(A)$ is
equal to
\begin{eqnarray*}
&&\sum_{i=0}^{k+1} \sum_{\sigma_i}
  (-1)^{s \sum\limits_{\ell \in \sigma_i^0} |\alpha_\ell|
  + \sum\limits_{j=1}^i \sum\limits_{\ell \in \sigma_i^0, \ell > \sigma_i(j)}
  |\alpha_{\sigma_i(j)}||\alpha_\ell|} \cdot  \\
&&\cdot \left ( \prod_{\ell \in \sigma_i^0} \frak
a_{-n_\ell}(\alpha_{\ell})
 \right )
 [[\cdots [\frak g, \frak a_{-n_{\sigma_i(1)}}(\alpha_{\sigma_i(1)})],
 \cdots], \frak a_{-n_{\sigma_i(i)}}(\alpha_{\sigma_i(i)})] |0\rangle
\end{eqnarray*}
where $0 \le i \le k+1$, $\sigma_i$ runs over all the maps $\{\,
1, \ldots, i \,\} \to \{\, 1, \ldots, b \,\}$ satisfying
$\sigma_i(1) < \cdots < \sigma_i(i)$, and $\sigma_i^0 = \{ \ell \,
| \, 1 \le \ell \le b, \ell \ne \sigma_i(1), \ldots, \sigma_i(i)
\}$.
\end{lemma}

\begin{definition} \label{partition}
Let $X$ be a smooth projective surface.
\begin{enumerate}
\item[{\rm (i)}] Let $\alpha \in H^*(X)$, and $\lambda = (\cdots
(-2)^{m_{-2}}(-1)^{m_{-1}} 1^{m_1}2^{m_2} \cdots)$ be a {\em
generalized partition} of the integer $n = \sum_i i m_i$ whose
part $i\in \Z$ has multiplicity $m_i$. Define $\ell(\lambda) =
\sum_i m_i$, $|\lambda| = \sum_i i m_i = n$, $\lambsq  = \sum_i
i^2 m_i$, $\lambda^! = \prod_i m_i!$, and
\begin{eqnarray*}
\mathfrak a_{\lambda}(\tau_*\alpha) = \left ( \prod_i (\mathfrak
a_i)^{m_i} \right ) (\tau_{\ell(\lambda)*}\alpha)
\end{eqnarray*}
where the product $\prod_i (\mathfrak a_i)^{m_i} $ is understood
to be $\cdots \mathfrak a_{-2}^{m_{-2}} \mathfrak a_{-1}^{m_{-1}}
 \mathfrak a_{1}^{m_{1}} \mathfrak a_{2}^{m_{2}}\cdots.$
Let $-\lambda$ be the generalized partition whose multiplicity of
$i \in \Z$ is $m_{-i}$.

\item[{\rm (ii)}] A generalized partition becomes a {\em partition}
in the usual sense if the multiplicity $m_ i = 0$ for every $i <
0$. A partition $\lambda$ of $n$ is denoted by $\lambda \vdash n$.

\item[{\rm (iii)}] We let ${\bf 1}_{-n}$ denote
$\mathfrak a_{-1}(1_X)^n/n!$ when $n \ge 0$ and $0$ when $n < 0$.
\end{enumerate}
\end{definition}

When $n \ge 0$, ${\bf 1}_{-n}\vac$ is the fundamental cohomology
class of the Hilbert scheme $\Xn$. The following theorem was one
of the main results proved in \cite{LQW4}.

\begin{theorem} \label{char_th}
Let $k \ge 0$ and $\alpha\in H^*(X)$. Then, $\mathfrak
G_k(\alpha)$ is equal to
\begin{eqnarray*}
& &- \sum_{\ell(\lambda) = k+2, |\lambda|=0}
   {1 \over \lambda^!} \mathfrak a_{\lambda}(\tau_{*}\alpha)
   + \sum_{\ell(\lambda) = k, |\lambda|=0}
   {\lambsq - 2 \over 24\lambda^!}
   \mathfrak a_{\lambda}(\tau_{*}(e\alpha))  \\
&+&\sum_{{\epsilon} \in \{\KX, \KX^2\}}  \,\,\,
   \sum_{\ell(\lambda) = k+2-|{\epsilon}|/2, |\lambda|=0}
   {g_\epsilon(\lambda) \over \lambda^!}
   \mathfrak a_{\lambda}(\tau_{*}(\epsilon\alpha))
\end{eqnarray*}
where all the numbers $g_\epsilon(\lambda)$ are independent of $X$
and $\alpha$.
\end{theorem}

\subsection{Ideals in $H^*(\Xn)$ for $X$ projective}
\par
${}$

\begin{lemma} \label{closed}
Let $\mathcal I$ be an ideal in the cohomology ring $H^*(X)$. Let
$\alpha \in \mathcal I$ and $k \ge 2$. Then, the pushforward
$\tau_{k*}\alpha$ can be written as $\sum_j \alpha_{j,1} \otimes
\cdots \otimes \alpha_{j, k} \in H^*(X^k)$ such that for each
fixed $j$, there exists some $\ell$ with $\alpha_{j,\ell} \in
\mathcal I$.
\end{lemma}
\begin{demo}{Proof}
First of all, note that if $\tau_{k*}\alpha = \sum_j \alpha_{j,1}
\otimes \alpha_{j,2} \otimes \cdots \otimes \alpha_{j, k}$, then
$\tau_{(k+1)*}\alpha = \sum_j (\tau_{2*}\alpha_{j,1}) \otimes
\alpha_{j,2} \otimes \cdots \otimes \alpha_{j, k}$. Therefore, by
induction, it suffices to prove the lemma for $k = 2$. Now the
case $k = 2$ follows from the observation that if we write
$\tau_{2*}(1_X)=\sum_i\alpha_i\otimes \beta_i$, then
$\tau_{2*}(\alpha) =\sum_i (\alpha\alpha_i) \otimes \beta_i$ with
$(\alpha\alpha_i) \in \mathcal I$.
\end{demo}

In view of the preceding lemma, we introduce the following
important definition.

\begin{definition} \label{pushforward_ideal}
Let $X$ be a smooth projective surface, and $\mathcal I$ be an
ideal in the cohomology ring $H^*(X)$. For $n \ge 1$, define $\In$
to be the subset of $\Hn$ consisting of the linear spans of
Heisenberg monomials of the form $\frak a_{-n_1}(\alpha_{1})
\cdots \frak a_{-n_b}(\alpha_{b}) \vac$ where $\alpha_i \in
\mathcal I$ for some $i$, and $n_1, \ldots, n_b$ are positive with
$\sum_\ell n_\ell = n$.
\end{definition}

\begin{lemma} \label{\Hn_I}
Let $\mathcal I$ be an ideal in the cohomology ring $H^*(X)$.
Then,
\begin{enumerate}
\item[{\rm (i)}] the linear subspace $\In$ in $\Hn$ is an ideal.

\item[{\rm (ii)}] $G_k(\alpha, n) \in \In$ if $\alpha \in \mathcal I$.
\end{enumerate}
\end{lemma}
\begin{demo}{Proof}
(i) Recall from \cite{LQW1} that the cohomology ring $\Hn$ is
generated by the cohomology classes $G_k(\alpha, n)$. So it
suffices to prove
%that $G_k(\alpha, n) \cdot \frak a_{-n_1}(\alpha_{1}) \cdots
%\frak a_{-n_b}(\alpha_{b}) \vac \in \In$, i.e.,
\begin{eqnarray}  \label{\Hn_I.1}
\mathfrak G_k(\alpha) \frak a_{-n_1}(\alpha_{1}) \cdots \frak
a_{-n_b}(\alpha_{b}) \vac \in \In
\end{eqnarray}
whenever $\alpha_1 \in \mathcal I$, and $n_1, \ldots, n_b$ are
positive with $\sum_\ell n_\ell = n$. For simplicity, put
$\mathfrak g = \mathfrak G_k(\alpha)$. Then, the operator
$\mathfrak g$ is of bi-degree $(0, 2k+|\alpha|)$, and satisfies
(\ref{nonsense1.1}) by Theorem~\ref{char_th}. Now we see from
Lemma~\ref{nonsense1} that $\mathfrak G_k(\alpha) \frak
a_{-n_1}(\alpha_{1}) \cdots \frak a_{-n_b}(\alpha_{b}) \vac$ is a
linear combination of expressions of the form
\begin{eqnarray} \label{\Hn_I.2}
\left ( \prod_{\ell \in \sigma_i^0} \frak
a_{-n_\ell}(\alpha_{\ell})
 \right )
 [[\cdots [\frak g, \frak a_{-n_{\sigma_i(1)}}(\alpha_{\sigma_i(1)})],
 \cdots], \frak a_{-n_{\sigma_i(i)}}(\alpha_{\sigma_i(i)})] |0\rangle
\end{eqnarray}
where $0 \le i \le k+1$, $\sigma_i$ maps the set $\{\, 1, \ldots,
i \,\}$ to the set $\{\, 1, \ldots, b \,\}$ with $\sigma_i(1) <
\cdots < \sigma_i(i)$, and $\sigma_i^0 = \{ \ell \, | \, 1 \le
\ell \le b, \ell \ne \sigma_i(1), \ldots, \sigma_i(i) \}$.

If $1 \in \sigma_i^0$, then (\ref{\Hn_I.2}) is contained in $\In$.
In the following, we assume $1 \not \in \sigma_i^0$. So $1 =
\sigma_i(1)$ since $\sigma_i(1) < \cdots < \sigma_i(i)$. By
Theorem~\ref{char_th} and Lemma~\ref{k_s}~(i), (\ref{\Hn_I.2}) is
a linear combination of expressions of the form
\begin{eqnarray} \label{\Hn_I.3}
\left ( \prod_{\ell \in \sigma_i^0} \frak
a_{-n_\ell}(\alpha_{\ell}) \right ) \mathfrak
a_{-\lambda}(\tau_{*}(\epsilon\alpha\alpha_{\sigma_i(1)} \cdots
\alpha_{\sigma_i(i)})) \vac
\end{eqnarray}
where ${\epsilon} \in \{1_X, e, \KX, \KX^2\}$, $\lambda \vdash
\sum_{j=1}^i n_{\sigma_i(j)}$ and $\ell(\lambda) =
k+2-|\epsilon|/2 - i$. By Lemma~\ref{closed}, the expression
(\ref{\Hn_I.3}) is contained in $\In$ since $\alpha_{\sigma_i(1)}
= \alpha_{1} \in \mathcal I$. It follows that (\ref{\Hn_I.2}) is
contained in $\In$. This proves (\ref{\Hn_I.1}).

(ii) Recall from the Corollary~4.8 in \cite{LQW4} that
$G_k(\alpha, n)$ is equal to
\begin{eqnarray}
& &\sum_{0 \le j \le k} \sum_{\lambda
   \vdash (j+1) \atop \ell(\lambda)=k-j+1}
   {(-1)^{|\lambda|-1} \over \lambda^! \cdot |\lambda|!}
   \cdot {\bf 1}_{-(n-j-1)} \mathfrak a_{-\lambda}(\tau_*\alpha)|0\rangle
   \nonumber \\
&+&\sum_{0 \le j \le k} \sum_{\lambda \vdash (j+1) \atop
\ell(\lambda)=k-j-1}
   {(-1)^{|\lambda|} \over \lambda^! \cdot |\lambda|!}
   \cdot {|\lambda| + \lambsq - 2 \over 24}
   \cdot {\bf 1}_{-(n-j-1)}
   \mathfrak a_{-\lambda}(\tau_*(e\alpha))|0\rangle  \nonumber \\
&+&\sum_{\epsilon \in \{\KX, \KX^2\} \atop 0 \le j \le k}
   \sum_{\lambda \vdash (j+1) \atop
\ell(\lambda)=k-j+1-|{\epsilon}|/2}
   {(-1)^{|\lambda|}g_{\epsilon}(\lambda+(1^{j+1}))
   \over \lambda^! \cdot |\lambda|!} \cdot {\bf 1}_{-(n-j-1)}
   \mathfrak a_{-\lambda}(\tau_*(\epsilon\alpha))|0\rangle
   \qquad \label{\Hn_I.4}
\end{eqnarray}
where $g_{\epsilon}$ is from Theorem~\ref{char_th}, and
$\lambda+(1^{j+1})$ is the partition obtained from $\lambda$ by
adding $(j+1)$ to the multiplicity of $1$. So $G_k(\alpha, n) \in
\In$ if $\alpha \in \mathcal I$.
\end{demo}

The following technical definition will be used throughout the
paper.

\begin{definition} \rm \label{universal}
Let $X$ be a smooth projective surface, $s \ge 1$, and $t_1,
\ldots, t_s \ge 1$. Fix $m_{i,j} \ge 0$ and $\beta_{i,j} \in
H^*(X)$ for $1 \le i \le s$ and $1 \le j \le t_i$. Then, a {\it
universal} linear combination of
$\displaystyle{\prod\limits_{j=1}^{t_i} G_{m_{i,j}}(\beta_{i,j},
n)}$, $1 \le i \le s$ is a linear combination of the form
$\displaystyle{\sum_{i=1}^s f_i \prod\limits_{j=1}^{t_i}
G_{m_{i,j}}(\beta_{i,j}, n)}$ where the coefficients $f_i$ are
independent of $X$ and $n$. A {\it universal} linear combination
of $\frak a_{-n_{i, 1}}(\beta_{i,1}) \cdots \frak a_{-n_{i,
t_i}}(\beta_{i,t_i})\vac$, $1 \le i \le s$ with $n_{i,j} \ge 1$
and $n_{i,1} + \ldots + n_{i,t_i} = n$ is defined in a similar
way.
\end{definition}

\begin{theorem} \label{\Hn_I_thm}
Let $X$ be a smooth projective surface, and $\mathcal I$ be an
ideal in the cohomology ring $H^*(X)$. If $\mathcal I$ is
homogeneous (i.e., $\mathcal I = \oplus_{i=0}^4 (\mathcal I \cap
H^i(X))$), then the ideal $\In$ is generated by the classes
$G_k(\alpha, n)$ with $\alpha \in \mathcal I$.
\end{theorem}
\begin{demo}{Proof}
Note that every Heisenberg monomial in the ideal $\In$ can be
written as $A = \displaystyle{{\bf 1}_{-(n - n_0)} \left (
\prod\limits_{i=1}^{s} \frak a_{-n_{i}}(\alpha_{i}) \right )
\vac}$ where $s \ge 1$, $n_1, \ldots, n_s \ge 1$, $n_0 =
\sum\limits_{i=1}^s n_i$, and $\alpha_\ell$ is contained in
$\mathcal I$ and homogeneous for some $\ell$. By
Lemma~\ref{\Hn_I}~(ii), it suffices to show that $A \in \In$ is a
universal finite linear combination of expressions of the form
\begin{eqnarray} \label{combination_of_G1}
\prod\limits_{j=1}^t G_{m_j}(\beta_j, n)
\end{eqnarray}
where $\sum\limits_{j=1}^t (m_j+1) \le n_0$, and $\beta_\ell \in
\mathcal I$ for some $\ell$.
%In fact, we will show that this linear combination is {\it universal}.

Use induction on $n_0$. When $n_0 = 1$, $s=n_1 = 1$. So $A= {\bf
1}_{-(n - 1)} \frak a_{-1}(\alpha_{1})|0\rangle = G_0(\alpha_1,
n)$ by (\ref{\Hn_I.4}). Hence the statement in the previous
paragraph holds for $n_0 = 1$.

Next assume $n_0 > 1$. Let $k_i = n_i -1$ for $1 \le i \le s$.
Then, $k_i \ge 0$ for every $i$. By the Theorem~4.1 and Lemma~5.1
in \cite{LQW3}, the cup product $\displaystyle{\prod_{i=1}^s
G_{k_i}(\alpha_i, n)}$ is equal to $\displaystyle{\left (
\prod_{i=1}^s {(-1)^{k_i} \over (k_i+1)!} \right )} A$ (defined to
be {\it the leading term}) plus a universal finite linear
combination of expressions ${\bf 1}_{-(n-\tilde n_0)} \left (
\prod\limits_{i=1}^{\tilde s} \frak a_{-\tilde n_{i}}(\tilde
\alpha_{i}) \right ) \vac$ where $\tilde \alpha_\ell \in \mathcal
I$ for some $\ell$, $\tilde s \ge 1$, $\tilde n_1, \ldots, \tilde
n_{\tilde s} \ge 1$, and $\sum_{i=1}^{\tilde s} \tilde n_{i} =
\tilde n_0 < \sum_{i=1}^s (k_i+1) = n_0$. By induction, $A$ is a
universal finite linear combination of expressions of the form
(\ref{combination_of_G1}).
\end{demo}

\begin{remark}
The assumption in Theorem~\ref{\Hn_I_thm} that the ideal $\mathcal
I \subset H^*(X)$ is homogeneous can be dropped when the surface
$X$ is simply connected.
\end{remark}

\subsection{Relation with the affine plane}
\par
${}$

In the following, we study the quotient ring $\Hn/\In$ when
$\mathcal I = \oplus_{\ell=1}^4 H^\ell(X)$. Note that $\Hn/\In$
has a linear basis consisting of Heisenberg monomials of the form
$\frak a_{-n_1}(1_X)^{r_1} \cdots \frak a_{-n_k}(1_X)^{r_k} \vac$
where $r_1, \ldots, r_k \ge 1$, and $0< n_1 < \ldots < n_k$ with
$\sum_{\ell=1}^k r_\ell n_\ell = n$. So we have an isomorphism of
vector spaces:
\begin{eqnarray} \label{iso}
\Phi: \bigoplus_{n \ge 0} \Hn/\In \to \mathbb Q[q_1, q_2, \ldots ]
\end{eqnarray}
where $\mathbb Q[q_1, q_2, \ldots ]$ is the polynomial ring in
countably infinitely many variables. Setting the degree of the
variable $q_i$ to be $i$, we see that $\Phi$ maps $\Hn/\In$ to the
homogeneous component of $\mathbb Q[q_1, q_2, \ldots ]$ of degree
$n$.

\begin{lemma} \label{char_class_modulo}
Let $\mathcal I = \oplus_{\ell=1}^4 H^\ell(X)$. Then, the quotient
ring $\Hn/\In$ is generated by the classes $G_k(1_X, n)$, $k = 0,
1, \ldots, n-1$. Moreover,
\begin{eqnarray} \label{char_class_modulo.1}
G_k(1_X, n) \equiv {(-1)^k \over (k+1)!} \cdot {\bf 1}_{-(n-k-1)}
\mathfrak a_{-(k+1)}(1_X)\vac \pmod {\In \, }.
\end{eqnarray}
\end{lemma}
\begin{demo}{Proof}
Since the cohomology ring $\Hn$ is generated by the classes
$G_k(\alpha, n)$ with $0 \le k < n$ and $\alpha \in H^*(X)$, the
first statement follows from Lemma~\ref{\Hn_I}~(ii). To prove
(\ref{char_class_modulo.1}), we note from (\ref{\Hn_I.4}) that the
leading term in $G_k(1_X, n)$ is ${(-1)^k \over (k+1)!} \cdot {\bf
1}_{-(n-k-1)} \mathfrak a_{-(k+1)}(\alpha)|0\rangle$ corresponding
to $j = k$, $\lambda \vdash (j+1)=k+1$ and
$\ell(\lambda)=k-j+1=1$. The other terms in $G_k(1_X, n)$ contain
$\tau_{i*}(\epsilon)$ with either $i \ge 2$ or $\epsilon = e, K,
K^2 \in \mathcal I$, and hence are contained in $\In$. This proves
(\ref{char_class_modulo.1}).
\end{demo}

\begin{theorem} \label{C^2}
Let $X$ be a smooth projective surface, and $\mathcal I =
\oplus_{\ell=1}^4 H^\ell(X)$. Then, the quotient $\Hn/\In$ is
isomorphic to the cohomology ring $H^*((\C^2)^{[n]})$.
\end{theorem}
\begin{demo}{Proof}
By Lemma~\ref{char_class_modulo}, the quotient ring $\Hn/\In$ is
generated by the classes $G_k(1_X, n)$, $k = 0, 1, \ldots, n-1$.
So by the Theorem~4.10 in \cite{Lehn}, it suffices to show that
via the isomorphism $\Phi$ in (\ref{iso}), the linear operator
$\mathfrak g_k$ on $\mathbb Q[q_1, q_2, \ldots ]$ induced by the
operator $\mathfrak G_k(1_X)$ on $\displaystyle{\bigoplus_{n \ge
0} \Hn}$ is given by
\begin{eqnarray} \label{C^2.1}
\mathfrak g_k = {(-1)^k \over (k+1)!} \sum_{n_1, \ldots, n_{k+1} >
0} \,\, q_{n_1+\ldots+n_{k+1}} \partial_{n_1} \cdots
\partial_{n_{k+1}}
\end{eqnarray}
where $\partial_i = {i{\partial \over \partial q_i}}$. Indeed, let
$A = \frak a_{-n_1}(1_X) \cdots \frak a_{-n_b}(1_X)\vac \in \Hn$
where $n_1, \ldots, n_b > 0$ with $\sum_\ell n_\ell = n$. By
Lemma~\ref{nonsense1}, $\mathfrak G_k(1_X)(A)$ is equal to
\begin{eqnarray} \label{C^2.3}
\sum_{i=0}^{k+1} \sum_{\sigma_i} \left ( \prod_{\ell \in
\sigma_i^0} \frak a_{-n_\ell}(1_X) \right ) [[\cdots [\mathfrak
G_k(1_X), \frak a_{-n_{\sigma_i(1)}}(1_X)], \cdots], \frak
a_{-n_{\sigma_i(i)}}(1_X)] |0\rangle
\end{eqnarray}
where for each fixed $i$, $\sigma_i$ runs over all the maps $\{\,
1, \ldots, i \,\} \to \{\, 1, \ldots, b \,\}$ satisfying
$\sigma_i(1) < \cdots < \sigma_i(i)$, and $\sigma_i^0 = \{ \ell \,
| \, 1 \le \ell \le b, \ell \ne \sigma_i(1), \ldots, \sigma_i(i)
\}$. Note from Lemma~\ref{k_s} that $[[\cdots [\mathfrak a_{t_1}
\cdots \mathfrak a_{t_r}(\tau_{r*}\alpha), \frak
a_{-n_{\sigma_i(1)}}(1_X)], \cdots], \frak
a_{-n_{\sigma_i(i)}}(1_X)] |0\rangle \in \In$ if $i \le r-2$ or
$\alpha \in \mathcal I$. Hence by (\ref{C^2.3}) and
Theorem~\ref{char_th}, $\mathfrak G_k(1_X)(A)$ equals
\begin{eqnarray*}
&&\qquad\qquad\qquad \sum_{\sigma_{k+1}}
  \left ( \prod_{\ell \in \sigma_{k+1}^0} \frak a_{-n_\ell}(1_X) \right )
   \\
&&\cdot \left [ \left [ \cdots \left [-\sum_{\ell(\lambda) = k+2,
|\lambda|=0} {1 \over \lambda^!} \mathfrak
a_{\lambda}(\tau_{*}1_X), \frak a_{-n_{\sigma_{k+1}(1)}}(1_X)
\right ], \cdots \right ], \frak a_{-n_{\sigma_{k+1}({k+1})}}(1_X)
\right ] \vac
\end{eqnarray*}
modulo $\In$. By Lemma~\ref{k_s}~(i) again, the preceding
expression is equal to
\begin{eqnarray*}
\sum_{\sigma_{k+1}} \left ( -\prod_{\ell =1}^{k+1}
(-n_{\sigma_{k+1}(\ell)}) \right ) \left ( \prod_{\ell \in
\sigma_{k+1}^0} \frak a_{-n_\ell}(1_X) \right ) \frak
a_{-n_{\sigma_{k+1}(1)}- \ldots -n_{\sigma_{k+1}({k+1})}}(1_X)
\vac.
\end{eqnarray*}

Therefore, we conclude that for the induced operator $\mathfrak
g_k$,
\begin{eqnarray*}
\mathfrak g_k(q_{n_1} \cdots q_{n_b}) &=&(-1)^k
\sum_{\sigma_{k+1}}
   \left ( \prod_{\ell =1}^{k+1} n_{\sigma_{k+1}(\ell)} \right )
   \left ( \prod_{\ell \in \sigma_{k+1}^0} q_{n_\ell} \right )
   q_{n_{\sigma_{k+1}(1)}+ \ldots +n_{\sigma_{k+1}({k+1})}}  \\
&=&{(-1)^k \over (k+1)!} \sum_{\tilde \sigma_{k+1}}
   \left ( \prod_{\ell =1}^{k+1} n_{\tilde \sigma_{k+1}(\ell)} \right )
   \left ( \prod_{\ell \in \tilde \sigma_{k+1}^0} q_{n_\ell} \right )
   q_{n_{\tilde \sigma_{k+1}(1)}+ \ldots +n_{\tilde \sigma_{k+1}({k+1})}}
\end{eqnarray*}
where $\tilde \sigma_{k+1}$ runs over all injective maps $\{\, 1,
\ldots, k+1 \,\} \to \{\, 1, \ldots, b \,\}$ and $\tilde
\sigma_{k+1}^0 = \{ \ell \, | \, 1 \le \ell \le b, \ell \ne \tilde
\sigma_{k+1}(1), \ldots, \tilde \sigma_{k+1}({k+1}) \}$ (so the
only difference between $\tilde \sigma_{k+1}$ and $\sigma_{k+1}$
is that we have dropped the condition $\sigma_{k+1}(1) < \cdots <
\sigma_{k+1}({k+1})$). Finally, the above formula for $\mathfrak
g_k(q_{n_1} \cdots q_{n_b})$ is equivalent to (\ref{C^2.1}).
\end{demo}
\section{Partial $n$-independence of structure constants
for $X$ projective} \label{sect_independence}

Given a finite set $S$ which is a disjoint union of subsets $S_0$
and $S_1$, we denote by ${\mathcal P}(S)$ the set of
partition-valued functions $\rho =(\rho(c))_{c \in S}$ on $S$ such
that for every $c \in S_1$, the partition $\rho(c)$ is required to
be {\it strict} in the sense that $\rho(c) =(1^{m_1(\rho(c))}
2^{m_2(\rho(c))} \ldots )$ with $m_r(\rho(c)) = 0$ or $1$ for all
$r \ge 1$.

Now let us take a linear basis $S= S_0 \cup S_1$ of $H^*(X)$ such
that $1_X, [x] \in S_0$, $S_0 \subset H^{\rm even}(X)$ and $S_1
\subset H^{\rm odd}(X)$. If we write $\rho =(\rho (c))_{c \in S}$
and $\rho(c) =(r^{m_r(\rho(c))})_{r \ge 1} =(1^{m_1(\rho(c))}
2^{m_2(\rho(c))} \ldots )$, then we put $\ell(\rho) = \sum_{c \in
S} \ell(\rho (c))
  = \sum_{c\in S, r\geq 1} m_r(\rho(c))$ and
\begin{eqnarray*}
\Vert \rho \Vert = \sum_{c \in S} |\rho (c)|
  =\sum_{c\in S, r\geq 1} r \cdot m_r(\rho(c)).
\end{eqnarray*}
Given $\rho \in \mathcal P(S)$ and $n \ge \| \rho \|+
\ell(\rho(1_X))$, we define $\tilde \rho \in \mathcal P(S)$ by
putting $m_r(\tilde{\rho}(c)) =m_r({\rho}(c))$ for $c \in S - \{
1_X \}$, $m_1(\tilde{\rho}(1_X)) =n -\| \rho \|- \ell(\rho(1_X))$,
and $m_r(\tilde{\rho}(1_X)) =m_{r-1} (\rho(1_X))$ for $r \geq 2$.
Note that $\| \tilde \rho \| = n$. We define ${\frak b}_\rho(n)
\in H^*(\Xn)$ by
\begin{eqnarray}
   {\frak b}_\rho(n)
&=&\frac1{\prod_{r \geq 1} (r^{m_r(\tilde{\rho}(1_X))}
   m_r(\tilde{\rho}(1_X))!)} \left ( \prod_{c \in S} \prod_{r \ge 1}
   {\mathfrak a}_{-r}(c)^{m_r(\tilde{\rho}(c))} \right ) \vac
   \qquad \label{b.1}   \\
&=&\frac{{\bf 1}_{-(n-\|\rho\|- \ell(\rho(1_X)))}}{\prod_{r \geq
2}
   (r^{m_r(\tilde{\rho}(1_X))} m_r(\tilde{\rho}(1_X))!)}
   \left ( \prod_{c \in S, r \ge 1 \atop c \ne 1_X \,\text{or } r > 1}
   {\mathfrak a}_{-r}(c)^{m_r(\tilde{\rho}(c))}
   \right ) \vac.  \qquad   \label{b.2}
\end{eqnarray}
where we fix the order of the elements $c \in S_1$ appearing in
the product $\displaystyle{\prod_{c \in S}}$ once and for all. For
$0 \le n < \| \rho \|+ \ell(\rho(1_X))$, we set ${\frak b}_\rho(n)
= 0$. This is consistent with (\ref{b.2}) and
Definition~\ref{partition}~(iii).
%Note from our definition of ${\bf 1}_{-(n-\|\rho\|- \ell(\rho(1_X)))}$
%that ${\frak b}_\rho(n)$ is equal to zero
%when $n < \| \rho \|+ \ell(\rho(1_X))$.
We remark that the only part in ${\frak b}_\rho(n)$ involving $n$
is the factor ${\bf 1}_{-(n-\|\rho\|- \ell(\rho(1_X)))}$ in
(\ref{b.2}) when $n \ge \| \rho \|+ \ell(\rho(1_X))$.

As a corollary to the theorem of Nakajima and Grojnowski
\cite{Gro, Na1, Na2}, $H^*(\Xn)$ has a linear basis consisting of
the classes:
\begin{eqnarray} \label{linear_basis}
{\frak b}_{\rho}(n), \qquad \rho \in \mathcal P(S) \,\,\, {\text
{\rm and }} \Vert \rho \Vert + \ell(\rho(1_X)) \le n.
\end{eqnarray}
Fix a positive integer $n$ and $\rho, \sigma \in {\mathcal P}(S)$
satisfying $\Vert \rho \Vert + \ell(\rho(1_X)) \le n$ and $\Vert
\sigma \Vert + \ell(\sigma(1_X)) \le n$. Then we can write the cup
product ${\frak b}_{\rho}(n) \cdot {\frak b}_{\sigma}(n)$ as
\begin{eqnarray} \label{str_constant}
{\frak b}_{\rho}(n) \cdot {\frak b}_{\sigma}(n) = \sum_{\nu \in
{\mathcal P}(S)} {a}_{\rho \sigma}^{\nu}(n) \,\, {\frak
b}_{\nu}(n)
\end{eqnarray}
where we have used ${a}_{\rho \sigma}^{\nu}(n) \in \mathbb Q$ to
denote the structure constants.

\begin{proposition}   \label{prop:poly}
Let $X$ be a smooth projective surface. The structure constants
${a}_{\rho \sigma}^{\nu}(n)$ of the cohomology ring $H^*(\Xn)$ are
polynomials in $n$ of degree at most
\begin{eqnarray} \label{prop:poly.1}
(\Vert \rho \Vert + \ell(\rho(1_X))) + (\Vert \sigma \Vert +
\ell(\sigma(1_X))) - (\Vert \nu \Vert + \ell(\nu(1_X))).
\end{eqnarray}
\end{proposition}
\begin{proof}
This is a consequence of the much more powerful Theorem~5.1 in
\cite{LQW3}. More explicitly, let $f(\rho) = \prod_{r \geq 2}
(r^{m_r(\tilde{\rho}(1_X))} m_r(\tilde{\rho}(1_X))!)$. Then,
$f(\rho)$ is independent of $n$. By the Theorem~5.1 in
\cite{LQW3}, the cup product $f(\rho) {\frak b}_\rho(n) \cdot
f(\sigma) {\frak b}_\sigma(n)$ is a linear combination of
expressions of the form:
\begin{eqnarray*}
{(n-\Vert \nu \Vert - \ell(\nu(1_X)))! \over (n-\Vert \nu \Vert -
\ell(\nu(1_X)) - i)!} \cdot f(\nu) {\frak b}_\nu(n),
\end{eqnarray*}
such that $ i \ge 0$, $(\Vert \nu \Vert + \ell(\nu(1_X)) + i) \le
(\Vert \rho \Vert + \ell(\rho(1_X))) + (\Vert \sigma \Vert +
\ell(\sigma(1_X)))$, and all the coefficients in this linear
combination are independent of $n$. It follows that all the
structure constants ${a}_{\rho \sigma}^{\nu}(n)$ are polynomials
in $n$ of degree at most (\ref{prop:poly.1}).
\end{proof}

To state our main result in this section, let $\mathcal I =
H^4(X)$ and $S_{\mathcal I} = \{ [x] \} \subset S$. Regard
$\mathcal P(S - S_{\mathcal I}) \subset \mathcal P(S)$. Then,
(\ref{str_constant}) implies that
\begin{eqnarray} \label{product_H^4(X).1}
{\frak b}_{\rho}(n) \cdot {\frak b}_{\sigma}(n) \equiv \,\,
\sum_{\nu \in \mathcal P(S - S_{\mathcal I})} {a}_{\rho
\sigma}^{\nu}(n) \,\, {\frak b}_{\nu}(n) \pmod {\In}.
\end{eqnarray}

\begin{theorem} \label{product_modulo_H^4}
Let $X$ be a smooth projective surface, and $\mathcal I = H^4(X)$.
Then, all the structure constants ${a}_{\rho \sigma}^{\nu}(n)$ in
(\ref{product_H^4(X).1}) are independent of $n$.
\end{theorem}

To prove this theorem, we need to establish four technical lemmas
first.

\begin{lemma} \label{pushforward_modulo_H^4}
If $k \ge 2$ and $\alpha$ is homogeneous, then $\tau_{k*}\alpha =
\sum_j \alpha_{j,1} \otimes \cdots \otimes \alpha_{j, k}$ where
for each $j$, either $|\alpha_{j,\ell}| = 4$ for some $\ell$, or
$0 < |\alpha_{j,\ell}| < 4$ for every $\ell$.
\end{lemma}
\begin{demo}{Proof}
Assume $|\alpha_{j,\ell}| < 4$ for every $\ell$. If
$|\alpha_{j,\ell}| = 0$ for some $\ell$, then $4(k-1)+ |\alpha| =
|\tau_{k*}\alpha| = \sum_{i=1}^k |\alpha_{j,i}| \le 3(k-1)$. So
$(k-1)+ |\alpha| \le 0$, contradicting to $k \ge 2$.
\end{demo}

\begin{lemma} \label{G_times_mononial}
Let $\mathcal I = H^4(X)$, $s \ge 0$, $n_1, \ldots, n_s > 0$,
$\tilde n = \sum_{\ell=1}^s n_\ell$, and $n \ge \tilde n$. Let
$\alpha, \alpha_1, \ldots, \alpha_s \in H^*(X)$ be homogeneous.
Assume $k + |\alpha| \ge 1$ and $n_\ell+ |\alpha_\ell| \ge 2$ for
every $\ell$. Then modulo $\In$, the cup product $G_k(\alpha, n)
\cdot \left ( {\bf 1}_{-(n- \tilde n)} \frak a_{-n_1}(\alpha_1)
\cdots \frak a_{-n_s}(\alpha_s)\vac \right )$ is a universal
linear combination of the basis (\ref{linear_basis}).
\end{lemma}
\begin{demo}{Proof}
By Lemma~\ref{\Hn_I}~(i) and (ii), the statement is trivial if one
of the classes $\alpha, \alpha_1, \ldots, \alpha_s \in H^*(X)$ is
contained in $\mathcal I$. So in the rest of the proof, we assume
that none of the classes $\alpha, \alpha_1, \ldots, \alpha_s$ is
contained in $\mathcal I$.

Our argument is similar to the proof of Lemma~\ref{\Hn_I}~(i). Put
$\mathfrak g = \mathfrak G_k(\alpha)$. Then, $B \, {\overset {\rm
def} =} \, G_k(\alpha, n) \cdot {\bf 1}_{-(n-\tilde n)} \frak
a_{-n_1}(\alpha_1) \cdots \frak a_{-n_s}(\alpha_s) \vac$ is equal
to
\begin{eqnarray*}
\displaystyle{1 \over (n-\tilde n)!} \mathfrak g \frak
a_{-1}(1_X)^{n-\tilde n} \frak a_{-n_1}(\alpha_1) \cdots \frak
a_{-n_s}(\alpha_s)\vac.
\end{eqnarray*}
By Lemma~\ref{nonsense1}, $B$ is a universal finite linear
combination of expressions of the form:
\begin{eqnarray} \label{G_times_mononial.1}
&{1 \over (n-\tilde n)!} {n-\tilde n \choose j} \frak
a_{-1}(1_X)^{n - \tilde n -j} \left ( \prod_{\ell \in \sigma_i^0}
\frak a_{-n_\ell}(\alpha_{\ell})
 \right )&    \nonumber \\
&\cdot [[\cdots [[[[\frak g, \underbrace{ \frak a_{-1}(1_X)],
\ldots], \frak a_{-1}(1_X)]}_{j \,\, \text{times}}, \frak
a_{-n_{\sigma_i(1)}}(\alpha_{\sigma_i(1)})],
 \cdots], \frak a_{-n_{\sigma_i(i)}}(\alpha_{\sigma_i(i)})]\vac \quad&
\end{eqnarray}
where $0 \le j \le n-\tilde n$, $0 \le i \le s$, $0 \le j+i \le
k+1$, $\sigma_i$ maps the set $\{\, 1, \ldots, i \,\}$ to the set
$\{\, 1, \ldots, s \,\}$ with $\sigma_i(1) < \cdots <
\sigma_i(i)$, and $\sigma_i^0 = \{ \ell \, | \, 1 \le \ell \le s,
\ell \ne \sigma_i(1), \ldots, \sigma_i(i) \}$. By
Theorem~\ref{char_th} and Lemma~\ref{k_s}~(i), we conclude that
\begin{eqnarray*}
[[\cdots [[[[\frak g, \underbrace{ \frak a_{-1}(1_X)], \ldots],
\frak a_{-1}(1_X)]}_{j \,\, \text{times}}, \frak
a_{-n_{\sigma_i(1)}}(\alpha_{\sigma_i(1)})],
 \cdots], \frak a_{-n_{\sigma_i(i)}}(\alpha_{\sigma_i(i)})]\vac
\end{eqnarray*}
is a universal finite linear combination of expressions $\mathfrak
a_{-\lambda}(\tau_{*}(\epsilon\alpha\alpha_{\sigma_i(1)} \cdots
\alpha_{\sigma_i(i)})) \vac$, where $\epsilon \in \{ 1_X, e, K,
K^2\}$, $\lambda \vdash j+n_{\sigma_i(1)}+\cdots+n_{\sigma_i(i)}$
and $\ell(\lambda) = k+2-|\epsilon|/2 - (j+i)$. So $B$ is a
universal finite linear combination of expressions
\begin{eqnarray} \label{G_times_mononial.2}
{\bf 1}_{-(n- \tilde n -j)} \left ( \prod_{\ell \in \sigma_i^0}
\frak a_{-n_\ell}(\alpha_{\ell}) \right ) \mathfrak
a_{-\lambda}(\tau_{*}(\epsilon\alpha\alpha_{\sigma_i(1)} \cdots
\alpha_{\sigma_i(i)})) \vac
\end{eqnarray}
where $\lambda \vdash j+n_{\sigma_i(1)}+\cdots+n_{\sigma_i(i)}$
and $\ell(\lambda) = k+2-|\epsilon|/2 - (j+i)$.

To prove our lemma, we see from (\ref{b.2}) that it suffices to
show that modulo ${\mathcal I}^{[\tilde n + j]}$, the part
$\displaystyle{\left ( \prod_{\ell \in \sigma_i^0} \frak
a_{-n_\ell}(\alpha_{\ell}) \right ) \mathfrak
a_{-\lambda}(\tau_{*}(\epsilon\alpha\alpha_{\sigma_i(1)} \cdots
\alpha_{\sigma_i(i)})) \vac}$ in (\ref{G_times_mononial.2}) does
not contain $\frak a_{-1}(1_X)$. Since $n_\ell + |\alpha_\ell| \ge
2$ for every $\ell$, this is equivalent to show that modulo
${\mathcal I}^{[\ell(\lambda)]}$, the part $\mathfrak
a_{-\lambda}(\tau_{*}(\epsilon\alpha \alpha_{\sigma_i(1)}\cdots
\alpha_{\sigma_i(i)}))\vac$ in (\ref{G_times_mononial.2}) does not
contain $\frak a_{-1}(1_X)$. By
Lemma~\ref{pushforward_modulo_H^4}, this is true if $\ell(\lambda)
\ge 2$. So let $\ell(\lambda) =1$. Then, we have
\begin{eqnarray*}
\mathfrak a_{-\lambda}(\tau_{*}(\epsilon\alpha\alpha_{\sigma_i(1)}
\cdots \alpha_{\sigma_i(i)})) \vac = \mathfrak
a_{-t}(\epsilon\alpha\alpha_{\sigma_i(1)} \cdots
\alpha_{\sigma_i(i)}) \vac
\end{eqnarray*}
where $t= |\lambda| = j+n_{\sigma_i(1)}+\cdots+n_{\sigma_i(i)}$
and $k+2-|\epsilon|/2 - (j+i) = 1$.

If $\mathfrak a_{-\lambda}(\tau_{*}(\epsilon\alpha
\alpha_{\sigma_i(1)}\cdots \alpha_{\sigma_i(i)}))\vac$ contains
$\frak a_{-1}(1_X)$, then we must have $t = 1$ and
$|\epsilon\alpha\alpha_{\sigma_i(1)} \cdots \alpha_{\sigma_i(i)}|
= 0$. So $j+n_{\sigma_i(1)}+ \cdots+ n_{\sigma_i(i)} =1$,
$\epsilon = 1_X$, and $|\alpha|=|\alpha_{\sigma_i(1)}| = \ldots =
|\alpha_{\sigma_i(i)}| = 0$. Thus, either $j=0$, $i=1$,
$n_{\sigma_1(1)}=1$ and $|\alpha_{\sigma_1(1)}| = 0$, or $j=1$ and
$i=0$. The first case contradicts to $n_{\sigma_1(1)}+
|\alpha_{\sigma_1(1)}| \ge 2$. In the second case, we see from
$k+2-|\epsilon|/2 - (j+i) = 1$ that $k=0$, contradicting to $k +
|\alpha| \ge 1$. So $\mathfrak
a_{-\lambda}(\tau_{*}(\epsilon\alpha \alpha_{\sigma_i(1)}\cdots
\alpha_{\sigma_i(i)}))\vac$ can not contain $\frak a_{-1}(1_X)$.
\end{demo}

\begin{lemma} \label{G_times_self}
Let $\mathcal I = H^4(X)$, $s \ge 1$, $k_1, \ldots, k_s \ge 0$,
$k_i + |\alpha_i| \ge 1$ for every $i$.
\begin{enumerate}
\item[{\rm (i)}] Modulo $\In$,
$\displaystyle{\prod_{i=1}^s G_{k_i}(\alpha_{i}, n)}$ is a
universal linear combination of (\ref{linear_basis}).

\item[{\rm (ii)}] When $n \ge n_0 \, {\overset {\rm def} =} \,
\sum\limits_{i=1}^s (k_i+1)$, the leading term in the cup product
$\displaystyle{\prod_{i=1}^s G_{k_i}(\alpha_{i}, n)}$ is equal to
$\displaystyle{\left ( \prod_{i=1}^s {(-1)^{k_i} \over (k_i+1)!}
\right )} \displaystyle{{\bf 1}_{-(n - n_0)} \left (
\prod\limits_{i=1}^{s} \frak a_{-(k_{i}+1)}(\alpha_{i}) \right )
\vac}$, which is equal to a universal multiple of ${\frak
b}_{\nu}(n)$ for some $\nu \in {\mathcal P}(S)$.
\end{enumerate}
\end{lemma}
\begin{demo}{Proof}
(i) Use induction on $s$. When $s = 1$, our statement follows from
Lemma~\ref{G_times_mononial} (take the integer $s$ there to be
$0$). Now, assume that the statement is true for $s-1$ with $s \ge
2$, i.e., $\displaystyle{\prod_{i=1}^{s-1} G_{k_i}(\alpha_{i},
n)}$ is a universal finite linear combination of the basis classes
${\frak b}_{\nu}(n)$, $\nu \in {\mathcal P}(S)$ and $\Vert \nu
\Vert + \ell(\rho(1_X)) \le n$. Note from (\ref{b.2}) that up to
some universal factor, every basis class ${\frak b}_{\nu}(n)$ is
of the form ${\bf 1}_{-(n-\tilde n)}\frak a_{-n_1}(\beta_1) \cdots
\frak a_{-n_s}(\beta_s)\vac$, where $n_1, \ldots, n_s > 0$,
$\tilde n = n_1 + \ldots+ n_s$, $n \ge \tilde n$, and $n_\ell +
|\beta_\ell| \ge 2$ for every $\ell$. So by
Lemma~\ref{G_times_mononial}, our statement for $s$ follows.

(ii) The statement about the leading term comes from the proof of
Theorem~\ref{\Hn_I_thm}, while the other statement about the
universal multiple follows from (\ref{b.2}).
\end{demo}

\begin{lemma} \label{tilde_a_product_G}
Let $\mathcal I = H^4(X)$. Then modulo $\In$, the basis element
${\frak b}_{\rho}(n)$ is a universal finite linear combination of
products of the form $\displaystyle{\prod_{j=1}^t
G_{m_j}(\beta_{j}, n)}$ where $m_j + |\beta_{j}| \ge 1$ for each
$j$, and $\sum_{j=1}^t (m_j+1) \le \Vert \rho \Vert +
\ell(\rho(1_X))$.
\end{lemma}
\begin{demo}{Proof}
By (\ref{b.2}), up to some universal factor, the basis class
${\frak b}_{\rho}(n)$ is of the form $\displaystyle{{\bf 1}_{-(n -
n_0)} \left ( \prod\limits_{i=1}^{s} \frak
a_{-(k_{i}+1)}(\alpha_{i}) \right ) \vac}$, where $n \ge n_0 \,
{\overset {\rm def} =} \, \Vert \rho \Vert + \ell(\rho(1_X)) =
\sum_{i=1}^s (k_i+1)$, and $k_i \ge 0$ and $k_i + |\alpha_i| \ge
1$ for every $i$. Now our result follows from an induction on
$n_0$ the same way as in the proof of Theorem~\ref{\Hn_I_thm}, and
Lemma~\ref{G_times_self}~(i) and (ii).
\end{demo}

\begin{demo}{Proof of Theorem~\ref{product_modulo_H^4}}
By Lemma~\ref{tilde_a_product_G}, ${\frak b}_\rho(n)$ is a
universal finite linear combination of expressions of the form
$\displaystyle{\prod\limits_{j=1}^{t_1} G_{m_{1,j}}(\beta_{1,j},
n)}$ where $m_{1,j}+|\beta_{1,j}| \ge 1$ for every $j$. Similarly,
${\frak b}_\sigma(n)$ is a universal finite linear combination of
expressions of the form $\displaystyle{\prod\limits_{j=1}^{t_2}
G_{m_{2,j}}(\beta_{2,j}, n)}$ where $m_{1,j}+|\beta_{1,j}| \ge 1$
for every $j$. Therefore, ${\frak b}_\rho(n) \cdot {\frak
b}_{\sigma}(n)$ is a universal finite linear combination of
expressions of the form $\displaystyle{\prod\limits_{j=1}^{t}
G_{m_{j}}(\beta_{j}, n)}$ where $m_{j}+|\beta_{j}| \ge 1$ for
every $j$. Now our result follows from
Lemma~\ref{G_times_self}~(i).
\end{demo}

We end this section with a lemma to be used in the next section.
For convenience, when $\ell(\rho)=1$, i.e., when the partition
$\rho(c)$ is a one-part partition $(r)$ for some $c\in S$ and is
empty for all the other elements in $S$, we will simply write
${\frak b}_{\rho}(n) ={\frak b}_{r,c}(n)$.

\begin{lemma} \label{1_partition}
Let $\mathcal I = H^4(X)$. Then modulo $\In$, the basis class
${\frak b}_{\rho}(n)$ is a universal finite linear combination of
products of the form $\displaystyle{\prod_{j=1}^t {\frak
b}_{r_j,c_j}(n)}$.
\end{lemma}
\begin{demo}{Proof}
Note from (\ref{b.2}) that ${\frak b}_{r,c}(n) = {\bf
1}_{-(n-r)}\frak a_{-r}(c)\vac$ if $c \ne 1_X$, while ${\frak
b}_{r,c}(n) = {\bf 1}_{-(n-r-1)}\frak a_{-(r+1)}(c)\vac$ if $c =
1_X$. So by Lemma~\ref{tilde_a_product_G}, it suffices to show
that if $c \in S$ and $k+|c| \ge 1$, then modulo $\In$, $G_k(c,
n)$ is a universal finite linear combination of products of the
form $\displaystyle{\prod_{j=1}^t \left ({\bf 1}_{-(n-r_j)}\frak
a_{-r_j}(c_j)\vac \right )}$ where $r_j \ge 1$, $c_j \in S$ and
$r_j + |c_j| \ge 2$.

Use induction on $k$. When $k=0$, we have $G_0(c, n) = {\bf
1}_{-(n-1)}\frak a_{-1}(c)\vac$ by formula (\ref{\Hn_I.4}). Next,
we prove that the statement in the preceding paragraph holds for
$k \ge 1$ by assuming that it is true for $0, \ldots, k-1$. By
Lemma~\ref{G_times_self}~(i) and (ii), modulo $\In$,
$\displaystyle{G_k(c, n) - {(-1)^k/(k+1)!} \cdot {\bf
1}_{-(n-k-1)} \frak a_{-(k+1)}(c)|0\rangle}$ is a universal linear
combination of the basis elements ${\frak b}_{\rho}(n)$ satisfying
$\Vert \rho \Vert + \ell(\rho(1_X)) < (k+1)$. By
Lemma~\ref{tilde_a_product_G}, modulo $\In$, each ${\frak
b}_{\rho}(n)$ is a universal finite linear combination of products
of the form $\displaystyle{\prod_{j=1}^t G_{k_j}(\gamma_{j}, n)}$
where $\gamma_{j} \in S$, $k_j + |\gamma_{j}| \ge 1$, and
$\sum\limits_{j=1}^t (k_j+1) \le \Vert \rho \Vert +
\ell(\rho(1_X))$. Note that $k_j < k$ for every $j$. So by
induction, we conclude that modulo $\In$, $G_k(c, n)$ is a
universal finite linear combination of products of the form
$\displaystyle{\prod_{j=1}^t \left ({\bf 1}_{-(n-r_j)}\frak
a_{-r_j}(c_j)\vac \right )}$ where $r_j \ge 1$, $c_j \in S$ and
$r_j + |c_j| \ge 2$.
\end{demo}

\section{The cohomology ring $H^*(\Xn)$ for $X$ quasi-projective}
\label{sect_quasi}

In this section, we will apply our results in previous sections to
smooth quasi-projective surfaces. Our terminology
``quasi-projective" means ``quasi-projective but not projective".
Recall from \cite{Na1} that for a smooth quasi-projective surface
$X$, the creation operators are modelled on the Borel-Moore
homology $H^{BM}_*(X)$, while the annihilation operators are
modelled on the ordinary homology $H_*(X)$. Then the Fock space of
the Heisenberg algebra is taken to be the direct sum over all $n$
of the Borel-Moore homology groups $H^{BM}_*(X^{[n]})$ \cite{Na1}.
Let $H_c^*(X)$ be the cohomology with compact support. Using the
the Poincar\' e dualities $\text{PD}: H^{4 - i}(X) \to
H^{BM}_i(X)$ and $\text{PD}: H_c^{4 - i}(X) \to H_i(X)$, we can
regard the creation operators $\frak a_{-n}(\alpha)$ with $n > 0$
as being modelled on $H^*(X)$ (i.e., $\alpha \in H^*(X)$), while
we can regard the annihilation operators $\frak a_{n}(\beta)$ with
$n > 0$ as being modelled on $H_c^*(X)$ (i.e., $\beta \in
H_c^*(X)$). Accordingly, with the help of the Poincar\' e duality
$\text{PD}: H^{4n - i}(X^{[n]}) \to H^{BM}_i(X^{[n]})$, from now
on we can take the Fock space to be the direct sum of the ordinary
cohomology groups $H^*(X^{[n]})$.

\subsection{The $n$-independence of the structure constants}
\par
${}^{}$

Let $X$ be a smooth quasi-projective surface embedded in a smooth
projective surface $\overline{X}$, and let $\iota: X \to
\overline{X}$ be the inclusion map. Then we have induced
embeddings $\iota_n: X^{[n]} \to \overline{X}^{[n]}$ for $n \ge
0$, and induced ring homomorphisms $\iota_n^*:
H^*(\overline{X}^{[n]}) \to H^*(X^{[n]})$. The maps $\iota^*$ and
$\iota_n^*$ are related by the following.

\begin{lemma} \label{iota_n^*}
Let notations be as in the preceding paragraph. Then,
\begin{eqnarray}  \label{iota_n^*.1}
\iota_n^*(\frak a_{-n_1}({\overline \alpha}_1) \cdots \frak
a_{-n_k}({\overline \alpha}_k)\vac ) = \frak
a_{-n_1}(\iota^*{\overline \alpha}_1) \cdots \frak
a_{-n_k}(\iota^*{\overline \alpha}_k)\vac
\end{eqnarray}
where $k \ge 0$, $n_1, \ldots, n_k > 0$, and $n_1+ \ldots + n_k
=n$.
\end{lemma}
\begin{demo}{Proof}
For $n \ge 0$, let $(\iota_n^*)^{BM}: H_i^{BM}(\overline{X}^{[n]})
\to H_i^{BM}(X^{[n]})$ be the natural map induced by the embedding
$\iota_n: X^{[n]} \to \overline{X}^{[n]}$. Then, it is well-known
that $(\iota_n^*)^{BM} \circ \text{PD} = \text{PD} \circ
(\iota_n^*)^{BM}$ (see \cite{Na2}). Combining with Nakajima's
constructions in \cite{Na1}, we obtain $\iota_{m_1+m_2}^* \frak
a_{-m_2}({\overline \alpha})(A) = \frak a_{-m_2}(\iota^*{\overline
\alpha}) \iota_{m_1}^*(A)$ for $m_1 > 0$, $m_2 > 0$ and $A \in
H^*(\overline{X}^{[m_1]})$. Applying this repeatedly, we obtain
(\ref{iota_n^*.1}).
\end{demo}

Next, assume that $\iota^*: H^*(\overline{X}) \to H^*(X)$ is
surjective. Let $\mathcal I = \ker(\iota^*)$. Fix a linear basis
$S$ of $H^*(\overline{X})$ as in Sect.~\ref{sect_independence}
such that $S$ contains a linear basis $S_{\mathcal I}$ of
$\mathcal I$ and $S_X {\overset {\rm def} =} \iota^*(S -
S_{\mathcal I})$ is a linear basis of $H^*(X)$. By
Lemma~\ref{iota_n^*}, $\ker(\iota_n^*) = {\mathcal I}^{[n]}$ which
is defined in Definition~\ref{pushforward_ideal}. Also, a linear
basis of $H^*(X^{[n]})$ is given by
\begin{eqnarray} \label{linear_basis_Y}
{\frak b}_{\rho_X}(n), \qquad \rho_X \in \mathcal P(S_X) \,\,\,
{\text {\rm and }} \Vert \rho_X \Vert + \ell(\rho_X(1_X)) \le n
\end{eqnarray}
where ${\frak b}_{\rho_X}(n)$ is defined in a similar way as in
(\ref{b.1}) and (\ref{b.2}). So for $\rho_X, \sigma_X \in
{\mathcal P}(S_X)$, we can write the cup product ${\frak
b}_{\rho_X}(n) \cdot {\frak b}_{\sigma_X}(n)$ as
\begin{eqnarray} \label{str_constant_Y}
{\frak b}_{\rho_X}(n) \cdot {\frak b}_{\sigma_X}(n) = \sum_{\nu_X
\in \mathcal P(S_X)} {a}_{\rho_X \, \sigma_X}^{\nu_X}(n) \,\,
{\frak b}_{\nu_X}(n)
\end{eqnarray}
where ${a}_{\rho_X \, \sigma_X}^{\nu_X}(n) \in \mathbb Q$ stands
for the structure constants.

\begin{definition} \label{S}
A smooth quasi-projective surface $X$ satisfies the {\it
S-property} if it can be embedded in a smooth projective surface
$\overline{X}$ such that the induced ring homomorphism
$H^*(\overline{X}) \to H^*(X)$ is surjective.
\end{definition}

\begin{theorem} \label{structure_constant_Y}
Let $X$ be a smooth quasi-projective surface satisfying the
S-property. Then all the structure constants ${a}_{\rho_X \,
\sigma_X}^{\nu_X}(n)$ in (\ref{str_constant_Y}) are independent of
$n$.
\end{theorem}
\begin{demo}{Proof}
Let notations be as above. Note that $\iota^*: (S - S_{\mathcal
I}) \to S_X$ is bijective. Define $\rho, \sigma \in {\mathcal
P}(S)$ by putting $m_r(\rho({\overline c})) =
m_r(\rho_X(\iota^*{\overline c}))$ and $m_r(\sigma({\overline c}))
= m_r(\sigma_X(\iota^*{\overline c}))$ when ${\overline c} \in (S
- S_{\mathcal I})$, and $m_r(\rho({\overline c})) = 0 =
m_r(\sigma({\overline c}))$ when ${\overline c} \in S_{\mathcal
I}$. By Theorem~\ref{product_modulo_H^4},
\begin{eqnarray*}
{\frak b}_{\rho}(n) \cdot {\frak b}_{\sigma}(n) \equiv \,\,
\sum_{\nu \in {\mathcal P}(S - \{[x]\})} {a}_{\rho
\sigma}^{\nu}(n) \,\, {\frak b}_{\nu}(n) \pmod
{{\overline{\mathcal I}}^{[n]}}
\end{eqnarray*}
where $\overline{\mathcal I} = H^4(\overline{X})$ and all the
${a}_{\rho \sigma}^{\nu}(n)$ are independent of $n$. Since
$\overline{\mathcal I} \subset \mathcal I$,
\begin{eqnarray} \label{structure_constant_Y.1}
{\frak b}_{\rho}(n) \cdot {\frak b}_{\sigma}(n) \equiv \,\,
\sum_{\nu \in {\mathcal P}(S - S_{\mathcal I})} {a}_{\rho
\sigma}^{\nu}(n) \,\, {\frak b}_{\nu}(n) \pmod {{\mathcal
I}^{[n]}}
\end{eqnarray}
where all the structure constants ${a}_{\rho \sigma}^{\nu}(n)$ are
independent of $n$. By Lemma~\ref{iota_n^*}, $\iota_n^*({\frak
b}_{\rho}(n)) = {\frak b}_{\rho_X}(n)$ and $\iota_n^*({\frak
b}_{\sigma}(n)) = {\frak b}_{\sigma_X}(n)$. Therefore, applying
$\iota_n^*$ to (\ref{structure_constant_Y.1}), we see that all the
structure constants ${a}_{\rho_X \, \sigma_X}^{\nu_X}(n)$ in
(\ref{str_constant_Y}) are independent of $n$.
\end{demo}

Thanks to Theorem~\ref{structure_constant_Y}, we will simply
denote the structure constants ${a}_{\rho_X \,
\sigma_X}^{\nu_X}(n)$ in (\ref{str_constant_Y}) by ${a}_{\rho_X \,
\sigma_X}^{\nu_X}$. Next, we study ring generators for the
cohomology ring $H^*(\Xn)$ when $X$ is a smooth quasi-projective
surface satisfying the S-property. For $\alpha \in H^*(X)$, define
$G_k(\alpha,n) = \iota_n^*G_k(\overline{\alpha},n)$ where
$\overline{\alpha} \in H^*(\overline{X})$ satisfies
$\iota^*\overline{\alpha} = \alpha$. This is independent of the
choice of $\overline{\alpha}$ by Theorem~\ref{\Hn_I_thm} and the
linearity of $G_k(\overline{\alpha},n)$ in $\overline{\alpha}$.

\begin{theorem}   \label{th:generator}
Let $X$ be a smooth quasi-projective surface embedded in a smooth
projective surface $\overline{X}$ such that the induced map
$H^*(\overline{X}) \to H^*(X)$ is surjective. Then, the cohomology
classes $G_k(\alpha,n)$, as $ 0 \leq k <n$ and $\alpha$ runs over
a linear basis of $H^*(X)$, form a set of ring generators of
$H^*(\Xn)$.
\end{theorem}
\begin{demo}{Proof}
Let notations be as above. By Lemma~~\ref{iota_n^*}, $H^*(\Xn)
\cong H^*(\overline{X}^{[n]})/{\mathcal I}^{[n]}$. Recall that the
classes $G_k(\overline{\alpha},n)$, $ 0 \leq k <n$ and
$\overline{\alpha} \in S \subset H^*(\overline{X})$, form a set of
ring generators of $H^*(\overline{X}^{[n]})$. It follows from
Theorem~\ref{\Hn_I_thm} that the classes
$\iota_n^*G_k(\overline{\alpha},n)$, $ 0 \leq k <n$ and
$\overline{\alpha} \in (S - S_{\mathcal I})$, form a set of ring
generators of $H^*(\Xn)$. In other words, the classes
$G_k(\alpha,n)$, $ 0 \leq k <n$ and $\alpha \in S_X = \iota^*(S -
S_{\mathcal I})$, form a set of ring generators of $H^*(\Xn)$.
Finally, note that $S_X$ is a linear basis of $H^*(X)$.
\end{demo}

\subsection{A universal ring}
\par
${}^{}$

Let $\mathfrak A = -\mathfrak a_1 ([x]_c)$ where $[x]_c \in
H_c^4(X)$ is the Poincar\'e dual of the homology class in $H_0(X)$
represented by a point in $X$.

\begin{lemma} \label{operator_A}
Let $X$ be a smooth quasi-projective surface. Then the linear map
$\mathfrak A: H^*(X^{[n+1]}) \rightarrow H^*(X^{[n]})$ is
surjective. In fact, it sends ${\frak b}_{\rho_X}(n+1)$ to ${\frak
b}_{\rho_X}(n)$.
\end{lemma}
\begin{proof}
Follows from the definition of the cohomology classes ${\frak
b}_{\rho_X}(n)$.
\end{proof}

\begin{corollary} \label{operator_A:cor}
Let $X$ be a smooth quasi-projective surface which satisfies the
S-property. Then $\mathfrak A :H^*(X^{[n+1]}) \rightarrow
H^*(X^{[n]})$ is a surjective ring homomorphism.
\end{corollary}
\begin{proof}
Follows immediately from Lemma~\ref{operator_A} and
Theorem~\ref{structure_constant_Y}.
\end{proof}

\begin{remark}
We conjecture that the surjective linear map $\mathfrak A
:H^*(X^{[n+1]}) \rightarrow H^*(X^{[n]})$ is a ring homomorphism
for an arbitrary smooth quasi-projective surface $X$. By
Lemma~\ref{operator_A}, this is equivalent to the Constant
Conjecture in \cite{Wa}.
\end{remark}

\begin{definition}
Let $X$ be a smooth quasi-projective surface satisfying the
S-property. We define the {\it FH ring $\mathcal G_X$ associated
to $X$} to be the ring with a linear basis given by the symbols
$\mathfrak b_{\rho_X}$, ${\rho_X} \in \mathcal P(S_X)$, with the
multiplication given by
\begin{eqnarray*}
{\frak b}_{\rho_X} \cdot {\frak b}_{\sigma_X} = \sum_{\nu_X \in
\mathcal P(S_X)} {a}_{\rho_X \, \sigma_X}^{\nu_X} \,\, {\frak
b}_{\nu_X}.
\end{eqnarray*}
where the structure constants ${a}_{\rho_X \, \sigma_X}^{\nu_X}$
come from (\ref{str_constant_Y}).
\end{definition}

For a smooth quasi-projective surface $X$ satisfying the
S-property, define a linear map $\mathfrak A_n: \mathcal G_X
\rightarrow H^*(\Xn)$ by sending ${\frak b}_{\rho_X}$ to ${\frak
b}_{\rho_X}(n)$. By the definition of $\mathcal G_X$ and
Theorem~\ref{structure_constant_Y}, $\mathfrak A_n$ is a
surjective ring homomorphism. We can illustrate
Theorem~\ref{structure_constant_Y} in terms of the following
commutative diagram of surjective ring homomorphisms:
\begin{eqnarray}  \label{comm_dia}
\CD
\ldots@>{=}>> \mathcal G_X @>{=}>> \mathcal G_X @>{=}>> \ldots\\
@VVV  @VV{\mathfrak A_{n+1}}V  @VV{\mathfrak A_n}V  @VVV\\
\ldots  @>{\mathfrak A}>> H^*(X^{[n+1]}) @>{\mathfrak A}>>
H^*(\Xn)
     @>{\mathfrak A}>> \ldots
\endCD
\end{eqnarray}

Next, we study the structure of the FH ring $\mathcal G_X$. For
fixed $r \ge 1$ and $c \in S_X$, we use ${\frak b}_{r,c}$ to
denote ${\frak b}_{\rho_X}$ where $\rho_X \in {\mathcal P}(S_X)$
is defined by taking $\rho_X(c)$ to be the one-part partition
$(r)$ and $\rho_X(c')$ to be empty for each $c' \ne c$.

\begin{theorem} \label{th_structure}
Let $X$ be a smooth quasi-projective surface which satisfies the
S-property. Then the FH ring $\mathcal G_X$ is isomorphic to the
tensor product $P \otimes E$, where $P$ is the polynomial algebra
generated by $\frak b_{r,c},\;c\in S_X \cap H^{\rm even}(X), r \ge
1$ and $E$ is the exterior algebra generated by ${\frak
b}_{r,c},\; c\in S_X \cap H^{\rm odd}(X), r \ge 1$.
\end{theorem}
\begin{proof}
Let notations be as above. First of all, we see from
Lemma~\ref{1_partition} that the FH ring $\mathcal G_X$ is
generated by the elements ${\frak b}_{r,c}$ with $r \ge 1$ and $c
\in S_X$. Also, note that $\mathcal G_X$ is super-commutative and
${\frak b}_{r,c}^2 = 0$ for $c\in S_X \cap H^{\rm odd}(X)$.

It remains to show that as $\rho =(r^{m_r(c)})_{c \in S_X, r \ge
1}$ runs over ${\mathcal P}(S_X)$, the monomials
$\displaystyle{\prod_{c\in S_X, r\ge 1}{\frak b}_{r, c}^{m_r(c)}}$
are linearly independent in $\mathcal G_X$. Assume
\begin{eqnarray*}
\displaystyle{\sum_{i \in I} d_i \prod_{c\in S_X, r\ge 1} {\frak
b}_{r, c}^{m^i_r(c)} =0}
\end{eqnarray*}
where $d_i \in \mathbb Q$ and $\rho_i=(r^{m^i_r(c)})_{c \in S_X, r
\ge 1}$ runs over a finite set $I$ of distinct elements in
${\mathcal P}(S_X)$. By (\ref{b.2}), we have $\frak b_{r, c}(n) =
{\bf 1}_{-(n-r-\delta_c)}\frak a_{-(r+\delta_c)}(c)\vac$ where
$\delta_c = 0$ if $c \ne 1_X$ and $\delta_c = 1$ if $c = 1_X$. So
we conclude that
\begin{eqnarray*}
 \quad \sum_{i \in I} d_i  \prod_{c\in S_X, r\ge 1}
 \left ( {\bf 1}_{-(n-r-\delta_c)}
 \frak a_{-(r+\delta_c)}(c)\vac \right )^{m^i_r(c)}
=\sum_{i \in I} d_i  \prod_{c\in S_X, r\ge 1}
 \frak b_{r, c}(n)^{m^i_r(c)} =0.
\end{eqnarray*}
Since $H^*(\Xn) \cong H^*(\overline{X}^{[n]})/{\mathcal I}^{[n]}$,
we see from Lemma~\ref{iota_n^*} that
\begin{eqnarray} \label{th_structure.1}
\quad \sum_{i \in I} d_i \cdot \prod_{\overline{c} \in S -
S_{\mathcal I}, r\ge 1} \left ( {\bf
1}_{-(n-r-\delta_{\iota^*\overline{c}})} \frak
a_{-(r+\delta_{\iota^*\overline{c}})}(\overline{c})\vac \right
)^{m^i_r(\iota^*\overline{c})} = w \in {\mathcal I}^{[n]}.
\end{eqnarray}

Take an integer $n$ large enough such that $n \ge n_i
\stackrel{\rm def}{=} \sum_{r,\overline{c} \in S - S_{\mathcal I}}
(r+\delta_{\iota^*\overline{c}}) m^i_{r}(\iota^*\overline{c})$ for
all $i \in I$. By the Theorem~5.1 in \cite{LQW3},
Eq.~(\ref{th_structure.1}) can be rewritten as
\begin{eqnarray} \label{th_structure.2}
\sum_{i \in I} d_i \left ( {\bf 1}_{-(n-n_i)} \left (
\prod_{\overline{c} \in S - S_{\mathcal I}, r\ge 1} (\frak
a_{-(r+\delta_{\iota^*\overline{c}})}
(\overline{c}))^{m^i_r(\iota^*\overline{c})} \right ) \vac + w_i
\right ) = w \in {\mathcal I}^{[n]}
\end{eqnarray}
where each $w_i$ is a universal finite linear combination of ${\bf
1}_{-(n - m)} \prod_{p=1}^N \frak a_{-m_{p}}(\g_{p}) \cdot
|0\rangle$ with $m = \sum_{p=1}^N m_{p} < n_i$ and $\g_p \in S$
for every $p$. Write $w$ and every $w_i$ as linear combinations of
the basis (\ref{linear_basis}). Since $(S - S_{\mathcal I}) \cap
S_{\mathcal I} = \emptyset$, we see from (\ref{th_structure.2})
that
\begin{eqnarray} \label{th_structure.3}
\sum_{i} d_i \cdot {\bf 1}_{-(n-n_i)} \left ( \prod_{\overline{c}
\in S - S_{\mathcal I}, r\ge 1} (\frak
a_{-(r+\delta_{\iota^*\overline{c}})}
(\overline{c}))^{m^i_r(\iota^*\overline{c})} \right ) \vac = 0
\end{eqnarray}
where $i$ satisfies $n_i = {\rm max}\{ n_j|j \in I\}$. Since the
partitions $\rho_i=(r^{m^i_r(c)})_{c \in S_X, r \ge 1} \in
{\mathcal P}(S_X)$ are distinct, all the coefficients $d_i$ in
(\ref{th_structure.3}) must be zero. By repeating the above
argument, we see that $d_i = 0$ for all $i \in I$.
\end{proof}

\begin{corollary} \label{cor_basis}
Let $X$ be a smooth quasi-projective surface which satisfies the
S-property. Then for $n \ge 1$, the cohomology ring $H^*(\Xn)$ is
generated by the classes $\frak b_{r,c}(n)$ where $1 \le r \le n$
and $c$ runs over a linear basis of $H^*(X)$.
\end{corollary}
\begin{proof}
Follows from Theorem~\ref{th_structure} and the observation that
the ring homomorphism $\mathfrak A_n$ in the commutative diagram
(\ref{comm_dia}) is surjective.
\end{proof}

We remark that Theorem~\ref{th:generator} provides a set of ring
generators for $H^*(\Xn)$. So Corollary~\ref{cor_basis} gives us a
second set of ring generators for $H^*(\Xn)$, which is parallel to
the set of ring generators for $H^*(\overline{X}^{[n]})$ found in
\cite{LQW2}.

\subsection{Examples of quasi-projective surfaces with the S-property}

\begin{example} \label{example_pt}
Let $\overline{X}$ be a projective surface and let $X$ be the
quasi-projective surface obtained from $X$ with a point removed.
It is easy to see that this smooth quasi-projective surface $X$
satisfies the S-property.
\end{example}

\begin{example} \label{example_quotient}
Let $\Gamma$ be a finite subgroup of $SL_2(\C)$. Let $X$ be the
minimal resolution of the simple singularity $\C^2/\Gamma$. It is
known that this smooth quasi-projective surface $X$ satisfies the
S-property. Moreover, $K_X = 0$.
\end{example}

\begin{example} \label{example_cotangent}
(The cotangent bundle of a smooth projective curve) Consider the
ruled surface $\overline{X} =\mathbb P(\mathcal O_C(-K_C) \oplus
\mathcal O_C)$ where $C$ is a smooth projective curve.
%We follow the notations in Hartshorne's book.
Let $\sigma$ be the section (to the projection $\overline{X} \to
C$) corresponding to the natural surjection $\mathcal O_C(-K_C)
\oplus \mathcal O_C \to \mathcal O_C(-K_C) \to 0$, and put $X =
\overline{X}-\sigma$. Then, $X$ is the total space of the
cotangent bundle of $C$, and $K_X = 0$.
%Note that $K_{\overline{X}} = -2 \sigma$.

We claim that $X$ satisfies the S-property. In fact, the following
general statement is true. Let $\overline{X}=\mathbb P(\mathcal
L_1\oplus \mathcal L_2)$ where $\mathcal L_1$ and $\mathcal L_2$
are two invertible sheaves over $C$. Let $\sigma$ (resp.
$\sigma'$) be the section of $\overline{X} \to C$ corresponding to
the natural surjection $\mathcal L_1\oplus \mathcal L_2 \to
\mathcal L_1 \to 0$ (resp. $\mathcal L_1\oplus \mathcal L_2 \to
\mathcal L_2 \to 0$). Put $X =\overline{X}-\sigma$. Then $X$
satisfies the S-property. To see this, let $X'
=\overline{X}-\sigma'$, and notice that $X$ and $X'$ are affine
bundles over $C$. Hence $X$ is homotopic to $C$, and $H^i(X) \cong
H^i(C)$ for every $i$. Therefore, to verify the S-property of $X$,
it remains to verify the surjectivities of the induced
homomorphisms $r_i: H^i(\overline{X}) \to  H^i(X)$ for $i = 1, 2$.
Consider the relative cohomology group $H^2(\overline{X}, X)$. By
the excision theorem, we obtain $H^2(\overline{X}, X) \cong
H^2(\overline{X}-\sigma', X-\sigma')= H^2(X', X'-\sigma)$. By the
Thom isomorphism, since $X'$ is an affine bundle over $C$ with
$\sigma$ being the zero section, we have $H^2(X', X'-\sigma) \cong
H^0(C) \cong \mathbb Q$. Hence $H^2(\overline{X}, X) \cong \mathbb
Q$. Now consider the exact sequence
\begin{eqnarray*}
H^1(\overline{X})\buildrel{r_1}\over \longrightarrow
H^1(X)\buildrel{\delta}\over\longrightarrow H^2(\overline{X},
X)\longrightarrow H^2(\overline{X})
\buildrel{r_2}\over\longrightarrow H^2(X).
\end{eqnarray*}
Since $H^2(X) \cong H^2(C) \cong \mathbb Q$ and $H^2(\overline{X})
\cong \mathbb Q\oplus \mathbb Q$, we conclude that the map
$\delta$ must be zero and the map $r_2$ must be surjective.
Therefore, $r_1$ is also surjective.
\end{example}

\section{Orbifold cohomology rings of symmetric products}
\label{sect_orb}

The orbifold cohomology ring of an orbifold was introduced in
\cite{CR}. Given an even-dimensional {\em compact} complex
manifold $X$, the orbifold cohomology rings $\orbsym$ of symmetric
products $X^n/S_n$ were studied in \cite{FG, QW, Uri} (also cf.
\cite{LS2, Ru1}). The axiomatic approach in \cite{QW} is
self-contained within the framework of the symmetric products,
while it is parallel to the study of the cohomology rings of
Hilbert schemes $\Xn$ when $X$ is a smooth projective surface. The
results obtained in the previous sections for the cohomology rings
of $\Xn$ when $X$ is a smooth quasi-projective surface are built
on the works \cite{LQW3, LQW4} (also cf. \cite{Lehn, LQW1}). As
observed in \cite{QW}, all the results in \cite{LQW3, LQW4} admit
exact counterparts in the orbifold cohomology rings of symmetric
products. This allows us to obtain readily the results on
$\orbsym$ when $X$ is an even-dimensional {\em non-compact}
complex manifold, which are the counterparts of those on the
cohomology rings of $\Xn$ when $X$ is a smooth quasi-projective
surface.

In this section, we will formulate and sketch these analogous
results for the orbifold cohomology rings $\orbsym$ when $X$ is
non-compact. We will not repeat the proofs for these analogous
results since the proofs are the same as in the Hilbert scheme
setup. For notational simplicity, we will assume below that $X$ is
a smooth quasi-projective surface which satisfies the S-property
(the assumption on the dimension of $X$ can be relaxed without
extra difficulty). We will use the results of \cite{QW} freely and
refer the reader to {\em loc. cit.} for details.

Let $\Xbar$ be a smooth projective surface which contains $X$ such
that the pullback map $\iota^*: H^*(\Xbar) \rightarrow H^*(X)$ is
surjective, where $\iota: X \rightarrow \Xbar$ is the inclusion
map. Recall from Sect.~5.1 in \cite{QW} that there exists a family
of ring products, denoted by $\circ_t$, on $\orbsymbar$ depending
on a rational (or complex) parameter $t$. When $t=1$, this
coincides with the original definition of orbifold product
\cite{CR} and when $t=-1$, this coincides with the modified
product in \cite{LS2, FG, Uri}. Put
\begin{eqnarray*}
\mathcal F_{\Xbar} = \bigoplus_{n=0}^\infty \orbsymbar.
\end{eqnarray*}
We can define the Heisenberg algebra acting irreducibly on
$\mathcal F_{\Xbar}$ with linear basis ${}^t \mathfrak p_n(\g)$,
where $n \in\Z$ and $\g \in H^*(\Xbar)$. The elements in this
linear basis satisfy the commutation relation (cf. Sections~3.2
and 5.2 in \cite{QW}):
\begin{eqnarray}  \label{eq:orbheis}
  [{}^t \mathfrak p_m(\alpha), {}^t \mathfrak p_n(\beta)]
   = t^{1/3} m \delta_{m,-n} \int_{\Xbar} (\alpha \beta) \cdot
   \text{Id}_{\mathcal F_{\Xbar}}.
\end{eqnarray}
When $t=-1$, this matches with the commutation relation of the
Heisenberg algebra associated to the Hilbert schemes (compare with
(\ref{eq:heis})).

The above definitions of the Fock space and of Heisenberg algebra
readily extends to $X$. Set $\mathcal F_X = \bigoplus_{n=0}^\infty
\orbsym.$ The creation operators ${}^t \mathfrak p_{-n}(\g)$,
where $n >0$ and $\g \in H^*(X)$, are defined in the same way as
before, while the annihilation operators ${}^t \mathfrak p_n(\g)$,
where $n >0$, are modelled on $\g \in H_c^*(X)$.

The definition in Sections 3.4 and 5.2 of \cite{QW} of the
cohomology classes $\eta_n(\alpha)$ and $O^k(\alpha,n)$ in
$\orbsym$ remain to be valid for $X$ quasi-projective without any
change. We further define $O_k(\alpha,n) =O^k(\alpha,n)/k!$. In
the same way, we define the operator ${}^t \mathfrak O^k(\alpha)$
(resp. ${}^t \mathfrak O_k(\alpha)$) in $\End (\mathcal F_X)$ to
be the orbifold product $\circ_t$ with the class $O^k(\alpha,n)$
(resp. with the class $O_k(\alpha,n)$) in $\orbsym$ for every $n
\ge 0$.

The inclusion map $\iota: X \rightarrow \Xbar$ induces an evident
surjective ring homomorphism $\jmath_n^* : \orbsymbar \rightarrow
\orbsym$ (note that $\jmath_1^* = \iota^*$ is surjective by
assumption). We have the analogue of Lemma~\ref{iota_n^*}, namely,
\begin{eqnarray}
\jmath_n^* (\mathfrak p_{-n_1}(\bar{\alpha}_1) \cdots
{}^t\mathfrak p_{-n_s}(\bar{\alpha}_s)\vac) = {}^t\mathfrak
p_{-n_1}(\iota^*\bar{\alpha}_1) \cdots {}^t\mathfrak
p_{-n_s}(\iota^*\bar{\alpha}_s)\vac
\end{eqnarray}
where $n_1, \ldots, n_s > 0$, and $n_1+ \ldots + n_s =n$, and
$\bar{\alpha}_i \in H^*(\Xbar)$. Here $\vac$ denotes $1 \in
H_{\text{\rm orb}}^*(pt)$ as usual. Given $\bar{\alpha} \in
H^*(\Xbar)$, we have by construction
\begin{eqnarray} \label{eq:simplepull}
\jmath_n^* (O_k(\bar{\alpha},n)) = O_k(\iota^*\bar{\alpha},n).
\end{eqnarray}

The following is the counterpart of Theorem~\ref{th:generator}.

\begin{proposition}  \label{prop:orbgenerator}
Let $X$ be a smooth quasi-projective surface embedded in a smooth
projective surface $\overline{X}$ such that the induced map
$H^*(\overline{X}) \to H^*(X)$ is surjective. Then, the cohomology
classes $O_k(\alpha,n)$, as $ 0 \leq k <n$ and $\alpha$ runs over
a linear basis of $H^*(X)$, form a set of ring generators of
$H^*_{\text{\rm orb}}(X^n/S_n)$.
\end{proposition}

Apparently, we can introduce a linear basis $\mathfrak
q_{\rho_X}(n)$ of $\orbsym$ in terms of the Heisenberg generators
${}^t\mathfrak p_n(\alpha)$, where $\rho_X \in \mathcal P(S_X)$
such that $\Vert \rho_X \Vert +\ell (\rho_X(1_X)) \leq n$, which
is the counterpart of the linear basis $\mathfrak b_{\rho_X}$ for
$H^*(\Xn)$. We write
\begin{eqnarray}  \label{str_cons_q}
\mathfrak q_{\rho_X}(n) \circ_t \mathfrak q_{\sigma_X}(n)
  =\sum_{\nu_X} p_{\rho_X \sigma_X}^{\nu_X}(n) \mathfrak q_{\nu_X}(n)
\end{eqnarray}
where $p_{\rho_X \sigma_X}^{\nu_X}(n)$ denotes the structure
constants for the orbifold product. The following proposition is
the counterpart of Theorem~\ref{structure_constant_Y}.

\begin{proposition}  \label{prop:const}
Let $X$ be a smooth quasi-projective surface which satisfies the
S-property. Then all the structure constants $p_{\rho_X \,
\sigma_X}^{\nu_X}(n)$ are independent of $n$.
\end{proposition}

\begin{remark}
Thanks to Proposition~\ref{prop:const}, we will denote $p_{\rho_X
\, \sigma_X}^{\nu_X}(n)$ simply by $p_{\rho_X \,
\sigma_X}^{\nu_X}$. It follows from Proposition~\ref{prop:const}
that we can introduce a universal ring $\mathcal U_X$ (referred to
as the FH ring again) with a linear basis given by the symbols
$\mathfrak q_{\rho_X}$ and multiplication given by $\mathfrak
q_{\rho_X} \circ \mathfrak q_{\sigma_X} =\sum_{\nu_X} p_{\rho_X
\sigma_X}^{\nu_X} \mathfrak q_{\nu_X}$. This FH ring $\mathcal
U_X$ governs the orbifold cohomology ring $(\orbsym, \circ_t)$ for
a fixed smooth quasi-projective surface with the S-property and
for every $n$. Similarly, we have a second set of ring generators
for $\orbsym$ which is the counterpart of
Corollary~\ref{cor_basis}.
\end{remark}

We introduce the linear isomorphisms $\Theta: \mathcal F_X
\rightarrow \fock$ and $\Theta_n: \orbsym \rightarrow H^*(\Xn)$ by
sending ${}^t\mathfrak p_{-n_1}(\alpha_1) \cdots {}^t\mathfrak
p_{-n_s}(\alpha_s)\vac$ to $\mathfrak a_{-n_1}(\alpha_1) \cdots
\mathfrak a_{-n_s}(\alpha_s)\vac$. Similarly, we define the linear
isomorphisms $\overline{\Theta}: \mathcal F_{\Xbar} \to \mathbb
H_{\Xbar}$ and $\overline{\Theta}_n: \orbsymbar \rightarrow
H^*(\Xbar^n/S_n)$. We have the following commutative diagram by
definitions:
\begin{eqnarray}  \label{dia:pullback}
\CD
\orbsymbar @>{\jmath_n^*}>> \orbsym  \\
  @VV{\overline{\Theta}_n}V  @VV{\Theta_n}V  \\
  H^*(\Xbar^{[n]}) @>{\iota_n^*}>>  H^*(\Xn)
\endCD
\end{eqnarray}

\begin{theorem}  \label{th:ringisom}
Let $X$ be a smooth quasi-projective surface with the S-property
and numerically trivial canonical class. Then the linear map
$\Theta_n: H^*_{\text{\rm orb}}(X^n/S_n) \to H^*(\Xn)$ is a ring
isomorphism, if we use the product $\circ_{-1}$ on $H^*_{\text{\rm
orb}}(X^n/S_n)$.
\end{theorem}
\begin{proof}
Set $t = -1$. However, for notational convenience, we will keep
writing ${}^t\mathfrak p$ instead of ${}^{-1}\mathfrak p$, etc.
The axiomatic approach in \cite{QW} shows that the operator
$\mathfrak O_k(\bar{\alpha}) \in \End (\mathcal F_{\Xbar})$, where
$\bar{\alpha} \in H^*(\Xbar)$, is equal to
\begin{eqnarray}  \label{eq:vertex}
& &- \sum_{\ell(\lambda) = k+2, |\lambda|=0}
   {1 \over \lambda^!} \,\,\,
   {}^t \mathfrak p_{\lambda}(\tau_{*}\bar{\alpha})
   + \sum_{\ell(\lambda) = k, |\lambda|=0}
   {\lambsq - 2 \over 24\lambda^!} \,\,\,
   {}^t \mathfrak p_{\lambda}(\tau_{*}(e_{\Xbar} \bar{\alpha}))
\end{eqnarray}
where $e_{\Xbar}$ is the Euler class of $\Xbar$. We remark that no
term in (\ref{eq:vertex}) involves the canonical class $K_{\Xbar}$
of $\Xbar$ in contrast to the formula for $\mathfrak
G_k(\bar{\alpha})$ in Theorem~\ref{char_th}.

Let ${\epsilon} \in \{K_{\Xbar}, (K_{\Xbar})^2\}, \ell(\lambda) =
k+2-|{\epsilon}|/2$ and $|\lambda|=0$. Let $\bar{\alpha}_i \in
H^*(\Xbar)$, $n_i> 0$ ($i =1, \ldots, s$), and $n_1+ \ldots + n_s
=n$. Applying the analogues of Lemma~\ref{k_s}~(i) to ${}^t
\mathfrak p$, we see that the expression ${}^t \mathfrak
p_{\lambda}(\tau_{*}(\epsilon \bar{\alpha})) \,\, {}^t\mathfrak
p_{-n_1}(\bar{\alpha}_1) \cdots {}^t\mathfrak
p_{-n_s}(\bar{\alpha}_s)\vac$ is a linear combination of
Heisenberg monomials of the form
\begin{eqnarray*}
{}^t\mathfrak p_{-n_1'}(\bar{\alpha}_1') \cdots {}^t\mathfrak
p_{-n_u'}(\bar{\alpha}_u') ({}^t\mathfrak p_{-n_{u+1}'} \cdots
{}^t\mathfrak p_{-n_{u+v}'})
(\tau_{v*}(K_{\Xbar}\bar{\alpha}'))\vac
\end{eqnarray*}
where $v > 0$, $n_1', \ldots, n_{u+v}' > 0$, and $n_1'+ \ldots +
n_{u+v}'=n$. From the proof of Lemma~\ref{closed},
$\tau_{v*}(K_{\Xbar}\bar{\alpha}')= \sum_j
(K_{\Xbar}\bar{\alpha}_{1, j}'') \otimes \bar{\alpha}_{2, j}''
\otimes \cdots \otimes \bar{\alpha}_{v, j}''$. Since
$\iota^*K_{\Xbar} = K_X = 0$ by assumption, $\iota^*(K_{\Xbar}
\bar{\alpha}_{1, j}'') = 0$. By Lemma~\ref{iota_n^*}, we conclude
that
\begin{eqnarray} \label{eq:zero}
\iota^*_n\overline{\Theta}_n \left ({}^t \mathfrak
p_{\lambda}(\tau_{*}(\epsilon \bar{\alpha})) \,\, {}^t\mathfrak
p_{-n_1}(\bar{\alpha}_1) \cdots {}^t\mathfrak
p_{-n_s}(\bar{\alpha}_s)\vac \right ) = 0.
\end{eqnarray}

Now consider a given Heisenberg monomial $A =\mathfrak
a_{-n_1}(\alpha_1) \cdots \mathfrak a_{-n_s}(\alpha_s)\vac \in
H^*(\Xn)$, where $\alpha_i \in H^*(X)$, $n_i> 0$ ($i =1, \ldots,
s$), and $n_1+ \ldots + n_s =n$. Fix $\bar{\alpha}_i \in
H^*(\Xbar)$ such that $\iota^*(\bar{\alpha}_i) =\alpha_i$. Put
$P={}^t\mathfrak p_{-n_1}(\alpha_1) \cdots {}^t\mathfrak
p_{-n_s}(\alpha_s)\vac$. Given $\alpha \in H^*(X)$, we choose
$\bar{\alpha} \in H^*(\Xbar)$ such that $\iota^*(\bar{\alpha})
=\alpha$. We have
\begin{eqnarray*}
\Theta_n ( O_k(\alpha,n) \circ_t P) \nonumber
  &=& \Theta_n \jmath_n^*( O_k(\bar{\alpha},n) \circ_t \,\,
  {}^t\mathfrak p_{-n_1}(\bar{\alpha}_1) \cdots {}^t\mathfrak
  p_{-n_s}(\bar{\alpha}_s)\vac)   \\
&=&\Theta_n \jmath_n^*( \mathfrak O_k(\bar{\alpha}) \,\,
  {}^t\mathfrak p_{-n_1}(\bar{\alpha}_1) \cdots {}^t\mathfrak
  p_{-n_s}(\bar{\alpha}_s)\vac)   \\
&=& \iota_n^* \overline{\Theta}_n ( \mathfrak O_k(\bar{\alpha})
  \,\, {}^t\mathfrak p_{-n_1}(\bar{\alpha}_1) \cdots
  {}^t\mathfrak p_{-n_s}(\bar{\alpha}_s)\vac)
\end{eqnarray*}
where we have used (\ref{dia:pullback}) and the fact that
$\jmath_n^*$ is a ring homomorphism. By (\ref{eq:vertex}),
(\ref{eq:zero}) and Theorem~\ref{char_th}, we get $\iota_n^*
\overline{\Theta}_n (\mathfrak O_k(\bar{\alpha}) \,\,
{}^t\mathfrak p_{-n_1}(\bar{\alpha}_1) \cdots {}^t\mathfrak
p_{-n_s}(\bar{\alpha}_s)\vac) = \iota_n^* ( \mathfrak
G_k(\bar{\alpha}) \mathfrak a_{-n_1}(\bar{\alpha}_1) \cdots
\mathfrak a_{-n_s}(\bar{\alpha}_s)\vac)$. Since $\iota_n^*$ is a
ring homomorphism, we obtain
\begin{eqnarray} \label{eq:isomring}
& &\Theta_n ( O_k(\alpha,n) \circ_t P) \nonumber
  = \iota_n^* ( \mathfrak G_k(\bar{\alpha})
  \mathfrak a_{-n_1}(\bar{\alpha}_1) \cdots \mathfrak
  a_{-n_s}(\bar{\alpha}_s)\vac) \nonumber \\
&=& \iota_n^* (G_k(\bar{\alpha}, n) \cdot
  \mathfrak a_{-n_1}(\bar{\alpha}_1) \cdots \mathfrak
  a_{-n_s}(\bar{\alpha}_s)\vac) = G_k(\alpha, n) \cdot A
\end{eqnarray}
noting that $\iota_n^*G_k(\bar{\alpha}, n) = G_k(\alpha, n)$ by
definition. Letting $P$ be the unit of the cohomology ring
$\orbsym$, we obtain from (\ref{eq:isomring}) that
$\Theta_n(O_k(\alpha,n)) =G_k(\alpha,n)$. Now our theorem follows
from (\ref{eq:isomring}), Theorem~\ref{th:generator} and
Proposition~\ref{prop:orbgenerator}.
\end{proof}

\begin{corollary}  \label{equal_str_cons}
With the same assumptions as in Theorem~\ref{th:ringisom}, then
the structure constant ${a}_{\rho_X \, \sigma_X}^{\nu_X}$ in
(\ref{str_constant_Y}) and the structure constant $p_{\rho_X \,
\sigma_X}^{\nu_X}$ in (\ref{str_cons_q}) are equal.
\end{corollary}
\begin{proof}
Follows immediately from (\ref{str_constant_Y}),
(\ref{str_cons_q}) and Theorem~\ref{th:ringisom}.
\end{proof}

\begin{remark} \label{rmk:ringisom}
The analogue to Theorem~\ref{th:ringisom} for smooth projective
surfaces with numerically trivial canonical class was established
in \cite{LS2, FG, Uri} and \cite{QW}. As pointed out in \cite{QW},
when the coefficient is taken to be $\C$ instead of $\mathbb Q$,
there exist explicit ring isomorphisms among the rings
$(H^*_{\text{orb}}(X^n/S_n; \C),\circ_t)$ for nonzero $t$. In
particular, when combined with Theorem~\ref{th:ringisom}, this
implies that the cohomology ring $H^*(\Xn; \C)$ is isomorphic to
the original orbifold cohomology ring $(H^*_{\text{orb}}(X^n/S_n;
\C),\circ_1)$ for smooth quasi-projective surfaces $X$ satisfying
the S-property and having numerically trivial canonical classes.
This further supports Ruan's conjecture in \cite{Ru1, Ru2}.
Finally, we notice that there are many examples of smooth
quasi-projective surfaces with the S-property and numerically
trivial canonical classes, including
Example~\ref{example_quotient} and
Example~\ref{example_cotangent}.
\end{remark}

\end{document}